\documentclass[12pt,a4paper]{article}

\usepackage[english]{babel}
\usepackage{amssymb}
\usepackage[latin1]{inputenc}
\usepackage{graphicx}
\usepackage{amsthm}
\usepackage{fullpage}

\title{Stability and instability of Navier boundary layers}
\date{September 2014}
\author{\textsc{Matthew Paddick} \\ IRMAR, Universit\'e de Rennes 1, France \\ \texttt{matthew.paddick@univ-rennes1.fr}}

\begin{document}

\renewcommand{\labelitemi}{$\bullet$}
\newtheorem{theo}{Theorem}[section]
\newtheorem{coro}[theo]{Corollary}
\newtheorem{lemma}[theo]{Lemma}
\newtheorem*{defi}{Definition}
\newtheorem{propo}[theo]{Proposition}
\newcommand{\preu}[1]{\textit{\underline{#1}}}
\newcommand{\eqref}[1]{(#1)}
\newcommand{\thref}[1]{Theorem \ref{#1}}
\newcommand{\lemref}[1]{Lemma \ref{#1}}
\newcommand{\propref}[1]{Proposition \ref{#1}}

\newcommand{\tilu}{\tilde{u}}
\newcommand{\tilv}{\tilde{v}}
\newcommand{\supp}{\mathrm{supp}}
\newcommand{\rplus}{\mathbb{R}^{+}}
\newcommand{\rrplus}{\mathbb{R}\times\rplus}
\newcommand{\bOmega}{\overline{\Omega}}
\newcommand{\dOmega}{\partial\Omega}
\renewcommand{\div}{\mathrm{div}~}
\newcommand{\rot}{\mathrm{rot}~}
\newcommand{\eps}{\varepsilon}
\newcommand{\ueps}{u^{\eps}}
\newcommand{\peps}{p^{\eps}}
\newcommand{\weps}{w^{\eps}}
\newcommand{\norme}[2]{\left\| #2 \right\| _{#1}}
\newcommand{\tilnorme}[2]{|| #2 ||_{#1}}
\newcommand{\tang}{_{\mathrm{tan}}}
\newcommand{\ust}{u_{sh}}
\newcommand{\bust}{\overline{\ust}}
\newcommand{\uap}{u^{ap}}
\newcommand{\pap}{p^{ap}}
\newcommand{\Rap}{R^{ap}}
\newcommand{\vap}{v^{ap}}
\newcommand{\deltaZ}{\delta_{0}}
\newcommand{\matdeux}[4]{\left(\begin{array}{cc} #1 & #2 \\ #3 & #4 \end{array} \right)}
\newcommand{\vecdeux}[2]{\left(\begin{array}{c} #1 \\ #2 \end{array} \right)}
\newcommand{\lapt}{{\cal L}_{t}}
\newcommand{\sdel}[1]{#1_{\delta}}
\newcommand{\rayl}[1]{(R(#1))}
\newcommand{\Teps}{T^{\eps}}
\newcommand{\uint}{u^{i}}
\newcommand{\ubl}{u^{b}}
\newcommand{\uD}{u^{D}}
\newcommand{\uN}{u^{N}}
\newcommand{\seta}[1]{#1_{\eta}}
\newcommand{\rotk}{\mathrm{rot}_k~}
\newcommand{\mathand}{\hspace{10pt} \mathrm{and} \hspace{10pt}}

\newcommand{\refrevu}[6]{\textsc{#1}, #2, \textit{#4} \textbf{#5} (#3), #6.}
\newcommand{\refbook}[4]{\textsc{#1}, \textit{#2} (#3), #4.}
\newcommand{\refbooks}[6]{\textsc{#1}, \textit{#2}, #3 \textbf{#4} (#5), #6.}

\maketitle

\begin{abstract}
We study the inviscid limit problem for the incompressible Navier-Stokes equation on a half-plane with a Navier boundary condition depending on the viscosity. On one hand, we prove the $L^{2}$ convergence of Leray solutions to the solution of the Euler equation. On the other hand, we show the nonlinear instability of some WKB expansions in the stronger $L^{\infty}$ and $\dot{H}^{s}$ ($s>1$) norms. These results are not contradictory, and in the periodic setting, we provide an example for which both phenomena occur simultaneously.
\newline

\end{abstract}

\section{Introduction}

\subsection{The inviscid limit problem}

We consider the two-dimensional incompressible Navier-Stokes equation on the half-plane $\Omega=\mathbb{R}\times]0,+\infty[$ with viscosity $\eps>0$:
$$ \eqref{NS(\eps)}: \hspace{10pt} \left\{ \begin{array}{rll} \partial_{t}\ueps + \ueps \cdot \nabla \ueps - \eps\Delta \ueps + \nabla \peps & = 0 & 
\mathrm{in}~\rplus\times\Omega \\ \div \ueps & = 0 & \mathrm{in}~\rplus\times\Omega \\ \ueps|_{t=0} & = \ueps_{0} & \mathrm{in}~\Omega . \end{array}\right. $$
$\ueps=(\ueps_{1},\ueps_{2})$ is the two-dimensional fluid speed and $\peps$ the kinematic pressure. We add two boundary conditions: first, the standard non-penetration condition for the component of $\ueps$ normal to the boundary
$$ \ueps\cdot \vec{n}=0 \hspace{10pt} \mathrm{on}~\rplus\times\partial\Omega , $$
where $\vec{n}=(0,-1)$, and a Navier condition that depends on the viscosity $\eps$ and that describes the tangential part of the fluid's speed on the boundary,
$$ \left(S\ueps\cdot \vec{n} + a^{\eps}\ueps\right)\tang=0 \hspace{10pt} \mathrm{on}~\rplus\times\partial\Omega , $$
where $a^{\eps}>0$, and $S\ueps=\frac{1}{2}(\nabla\ueps + ^{t}\nabla\ueps)$.

Unlike the homogeneous Dirichlet (no-slip) boundary condition, $u|_{y=0}=0$, the Navier (slip) condition, introduced in the first half of the XIX\textsuperscript{th} century by H.Navier himself (\cite{Nh}), allows the fluid to slide along the boundary. Physically, this is more accurate than the no-slip condition when one takes into account interaction at the boundary.
 The slip condition is therefore used to model, for example, blood flow in capillary vessels that are a few microns wide, and where molecular exchanges with the neighbouring cells take place (see \cite{PRD}).
 One also encounters a Navier condition when homogenising the no-slip condition on rough or porous walls (see \cite{JM} and \cite{GVM}). Mathematically, it is shown in \cite{MSR} that the Navier condition can be derived by taking the limit when the mean-free path goes to zero of renormalised solutions of the Boltzmann equation with a Maxwell reflection boundary condition.
\newline

This chapter will deal with the inviscid limit problem for $\eqref{NS(\eps)}$, \textit{i.e.} we study the behaviour of a family of solutions of $\eqref{NS(\eps)}$ with a Navier boundary condition, $(\ueps)_{0<\eps\leq\eps_{0}}$, relative to $v$, a solution of the incompressible Euler equation
$$ \eqref{E}: \hspace{10pt} \left\{ \begin{array}{rcl} \partial_{t}v + v\cdot\nabla  v + \nabla q & = & 0 \\ \div v & = & 0 \\ v|_{t=0} & = & v_{0} \\ v_{2}|_{y=0} & = & 0 \end{array}\right. $$
when the viscosity $\eps$ goes to 0.

The inviscid limit problem is well understood in the whole space or with periodic boundary conditions (see \cite{MB}). However, in domains with boundaries, a lack of compactness appears due to the presence of boundary layers, making the problem considerably more difficult. In recent years, much progress has been made in the case of non-characteristic boundaries.
 For incompressible fluids, these appear in the case of injection or succion boundary conditions (\cite{TW}). For related results for compressible models and general hyperbolic-parabolic systems, we refer to \cite{MZ}.
 In this situation the boundary layer is of size $\eps$ and has an amplitude ${\cal O}(1)$, and sharp stability and instability conditions can be shown (see for example \cite{DG} and \cite{Rf} for the study of the stability or instability of the Ekman layer in rotating fluids).

On the half-plane, the inviscid limit problem remains mainly unsolved in the case of Navier-Stokes equations with the homogeneous Dirichlet boundary condition. Formally, one expects the boundary layer to be of size $\sqrt{\eps}$ and hence to write a solution $\uD$ as
\begin{equation}
\uD(t,x,y)=\uint(t,x,y)+\ubl(t,x,\eps^{-1/2}y) , \label{anzD} \end{equation}
in which $\uint$ solves $\eqref{E}$ and $\ubl$ is a corrective term. A major difficulty resides in the existence of such an expansion, as, when one formally rescales the Navier-Stokes equation with $(u_{1},u_{2})=(\tilu,\sqrt{\eps}\tilv)$ in the variables $(t,x,y'=\eps^{-1/2}y)$, one gets that the behaviour of the solution near the boundary is governed by the Prandtl equation:
$$ \left\{ \begin{array}{rcl} \partial_{t}\tilu + \tilu\partial_{x}\tilu + \tilv\partial_{y'}\tilu-\partial^{2}_{y'y'}\tilu & = & (\partial_{t}\uint+\uint\cdot\nabla \uint)_{1}|_{y=0} \hspace{5pt} \mathrm{in}~\Omega \\
\partial_{x}\tilu + \partial_{y'}\tilv & = & 0 \hspace{5pt} \mathrm{in}~\Omega \\ (\tilu,\tilv)|_{y=0} & = & 0 \\ \lim_{y\rightarrow+\infty} \tilu & = & \uint|_{y=0} . \end{array} \right. $$
While the Prandtl system is well-posed with monotonous initial conditions in the strip $[0,L]\times\rplus$ (a fact known since the 60s, \cite{Oo}), it has recently been proved by D.G\'erard-Varet and E.Dormy in \cite{GVD} that the Prandtl equation is ill-posed in Sobolev spaces on $\mathbb{T}\times\rplus$.
 This shows that the ansatz (\ref{anzD}) is only formal in general, and E.Grenier showed in \cite{Ge} that, even when such an ansatz can be justified, \textit{e.g.} in analytic framework (see \cite{SC}), it may not be valid in $H^{1}$, and instability occurs on the derivatives in $L^{\infty}$. The instability phenomenon is linked to the linear instability of shear flows for the Euler equation.
\newline

For the Navier-Stokes equation with a Navier boundary condition, when $a^\eps$ is a fixed number independent of $\eps$, the inviscid limit problem is solved in $L^2$ framework, see \cite{Bc}, \cite{CMR}, \cite{Kjp}, \cite{LFNLP}, \cite{IP}. Moreover, asympotic expansions of the form
$$ \uN(t,x,y)=\uint(t,x,y)+ \sqrt{\eps}\ubl(t,x,\eps^{-1/2}y) , $$
with the amplitude of the boundary layer much smaller than in (\ref{anzD}), are rigorously justified in \cite{IS} and uniform conormal estimates in agreement with this behaviour are obtained in \cite{MR} and \cite{GK}. We also refer to \cite{XX}, \cite{BC2} and \cite{BC1} for the study of special cases in 3D.
   
Here, we shall be interested in the case where the Navier condition is as above with $a^\eps=a\eps^{-1/2}$, $a>0$ constant. Denoting the two components of the fluid's speed by $\ueps_{1}$ and $\ueps_{2}$, the boundary conditions translate as the following:
$$ \begin{array}{rrl} \eqref{NP}: & \ueps_{2}(t,x,0) & = 0 \\
(NC(a,\eps)): & \frac{1}{2}\partial_{y}\ueps_{1}(t,x,0) & = \frac{a}{\sqrt\eps} \ueps_{1}(t,x,0) . \end{array} $$
Our motivation for this is to understand the transition between the unstable Dirichlet case and the stable Navier case studied in \cite{IS}. Note that for $\eps$ small, our slip condition seems to approximate the no-slip Dirichlet condition, and we expect the formal asymptotic behaviour of the solution to replicate the ansatz of the Dirichlet case, \textit{i.e.}
$$ \uN(t,x,y)=\uint(t,x,y)+\ubl(t,x,\eps^{-1/2}y) , $$
hence we should observe some instability.
\newline

\textbf{Vocabulary and summary of the results.} As of now, the term ``solution of the Navier-Stokes equation $(NS(a,\eps))$'' will designate a solution of $\eqref{NS(\eps)}$ satisfying the boundary conditions $\eqref{NP}$ and $(NC(a,\eps))$.
 We obtain two results: firstly, convergence in $L^{2}$ of every sequence of Leray solutions of $(NS(a,\eps))$ (the full definition of which is given at the beginning of section 2) to a smooth solution of $\eqref{E}$ when $\eps$ goes to zero and the initial conditions are $H^{s}$ for some $s>2$,
 and secondly, we prove that there exist boundary layer profiles such that the WKB expansion, which is also of the form (\ref{anzD}) in our case, is unstable in $L^{\infty}$.

\subsection{The $L^{2}$ stability result}

We show that $\ueps$ converges to $v$ in $L^{2}(\Omega)$, extending the result obtained by D.Iftimie and G.Planas in \cite{IP} for a constant Navier boundary condition, \textit{i.e.} $(Su\cdot\vec{n}+a u)\tang=0$ with $a>0$ not depending on the viscosity.

\begin{theo} Let $\ueps$ be the Leray solution of the Navier-Stokes equation $(NS(a,\eps))$, and $v$ be the solution of the Euler equation $\eqref{E}$. We assume that the initial condition $v_{0}$ is in $H^{s}(\Omega)$ for some $s>2$. If $\ueps_{0}$ converges to $v_{0}$ in $L^{2}(\Omega)$ as $\eps$ goes to 0, then for every $T>0$,
$$ \sup_{t\in[0,T]} \norme{L^{2}(\Omega)}{\ueps(t)-v(t)} \stackrel{\eps\rightarrow 0}{\longrightarrow} 0 . $$ \label{theo1} \end{theo}
The proof, based on classical energy estimates, is as in \cite{IP}. This method also allows us to expose a convergence rate. Providing $\norme{L^2}{\ueps(0)-v(0)} = {\cal O}(\eps^{1/4})$, we get that $\norme{L^2}{\ueps(t)-v(t)}$ is also ${\cal O}(\eps^{1/4})$.
\newline

We single out $a^\eps=a\eps^{-1/2}$ because this is the case that will provide the nonlinear instability in the next paragraph, but we can extend the proof of this theorem to a whole family of Navier boundary conditions, of the type
$$ (NC(a,\eps,\beta)): \hspace{10pt} \partial_{y}\ueps_{1}(t,x,0) = \frac{2a}{\eps^{\beta}}\ueps_{1}(t,x,0) , $$
whatever the sign of $a$ (although the physical meaning of $a<0$ is not well understood). Precisely, we have convergence in the following cases:

\begin{theo} 
\label{theostab} Let $\ueps$ be the Leray solution of $(NS(a,\eps,\beta))$, the Navier-Stokes system in which the boundary condition $(NC(a,\eps))$ has been replaced by $(NC(a,\eps,\beta))$. Let $v$ and $v_{0}$ be as in \thref{theo1}, and $\ueps_{0}$ converge to $v_{0}$ in $L^{2}$ as $\eps\rightarrow0$. Then we have convergence of $\ueps$ to $v$ in the same sense as \thref{theo1} if:

- either $a>0$ and $\beta<1$,

- or $a<0$ and $\beta\leq 1/2$. \\
In both cases, the convergence rate is ${\cal O}(\eps^{(1-\beta)/2})$. \end{theo}

This extends a recent result by X-P.Wang, Y-G.Wang and Z.Xin in \cite{WWX}, in which they proved convergence for $\beta<1/2$, regardless of the sign of $a$. Notably, we add the case $\beta=1/2$, which seems critical when $a<0$. Also, when $a>0$, we can go further than $\beta=1/2$, and all the way up to $\beta<1$. The convergence rate is in agreement with other works (\cite{IP}, \cite{GK}).
\newline

We also point out that these results can be adapted to higher dimensions, for smoother limit data ($v_0\in H^s$ for $s>d/2+1$ in $d$ dimensions), and for as long as a strong solution to the Euler equation exists. We will explain in the proof how to get the starting energy estimate in higher dimensions.

\subsection{The nonlinear instability result}

In \cite{Ge}, E.Grenier proved that linear instability implies nonlinear instability for the Euler equation, and for some viscous boundary layers (\cite{DG} with B.Desjardins). It is this argument that we shall use to get our instability result, and we shall use a shear flow that is linearly unstable for the Euler equation as a starting point:

\begin{defi} A \underline{shear flow} is a 2D smooth vector field written as $\ust=(\ust(y),0)$ (the vector field and its first component are indifferently called $\ust$). Note that a shear flow, associated with a constant pressure, is automatically a stationary solution of the Euler equation.
 The shear flow $\ust$ is \underline{linearly unstable} if there exist $k\in\mathbb{R}$, $\lambda\in\mathbb{C}$ with $\Re(\lambda)>0$, and a function $\Psi\in H^{1}_{0}(\rplus)$ such that
$$ u(t,x,y) = \left(\begin{array}{c} e^{\lambda t}e^{ikx}\Psi'(y) \\ -ik e^{\lambda t} e^{ikx} \Psi(y) \end{array} \right) $$
is a solution of the Euler equation linearised around $\ust$,
$$ \eqref{EL}: \hspace{10pt} \left\{ \begin{array}{rcl} \partial_{t}u+\ust\cdot\nabla u+u\cdot\nabla\ust +\nabla p & = & 0 \\ \div u & = & 0 \end{array} \right. \hspace{10pt} \mathrm{in}~\rplus\times\Omega , $$
with $\eqref{NP}$ condition on the boundary. \end{defi}

\begin{theo} Let $\ust$ be a linearly unstable shear flow for the Euler equation. Setting $Y=\eps^{-1/2}y$, we generate a time-dependent boundary layer $\bust$ as solution of the heat equation
$$ \left\{ \begin{array}{rcl} \partial_{t}\bust(t,Y) - \partial^{2}_{YY}\bust(t,Y) & = & 0 \\ \bust(0,Y) & = & \ust(Y) \\ \frac{1}{2}\partial_{Y}\bust(t,0) & = & a\bust(t,0) , \end{array} \right. $$
the boundary condition being dictated by $(NC(a,\eps))$, $a>0$; $\bust(t,y/\sqrt{\eps})$ therefore solves the Navier-Stokes system (with constant pressure). Then, for any $n\in\mathbb{N}^*$, there exist $\deltaZ$ and $\eps_{0}>0$ such that for any $0<\eps<\eps_{0}$, there exists a solution $\ueps$ of $(NS(a,\eps))$ with initial data $u_{0}^\eps$ such that
$$ \norme{H^{s}(\Omega)}{\ueps_{0}-\ust\left(\frac{y}{\sqrt{\eps}}\right)}\leq C\eps^{n} $$
for some $s>0$, that satisfies the following: at a time $\Teps\sim n\ln(\eps^{-1})\sqrt{\eps}$,
\begin{equation} \norme{L^{\infty}(\Omega)}{\ueps(\Teps,x,y)-\bust\left(\Teps,\frac{y}{\sqrt{\eps}}\right)} \geq \deltaZ . \label{inst} \end{equation}
Moreover, for $\sigma>1$,
$$ \norme{\dot{H}^{\sigma}}{\ueps(\Teps,x,y)-\bust\left(\Teps,\frac{y}{\sqrt{\eps}}\right)} \stackrel{\eps\rightarrow 0}{\longrightarrow} +\infty . $$
\label{theo3} \end{theo}

Note that the solutions to the Navier-Stokes equation that we consider in this theorem are $L^{2}(\Omega)$ perturbations of $\bust$, which is itself not necessarily in $L^2(\Omega)$.
 However, we point out that, like the equivalent statement for Dirichlet boundary conditions in \cite{Ge}, a similar result holds in the domain $\Omega'= \mathbb{T} \times  ]0,+\infty[$, the proof being slightly easier, and in this case, we can have shear flows that are square-integrable, and then we expect both our main results to occur simultaneously:
 the instability of WKB expansions is not an obstruction to $L^2$-convergence. In the final section, we will provide an explicit example of a profile which will lead to this situation.
\newline

What \thref{theo3} provides is a nonlinear instability result in $L^\infty$ in the usual sense for WKB expansions of the form (\ref{anzD}): we show that when the internal layer $u^i$ is linearly unstable for the Euler equation, then the WKB expansion is unstable, in the sense that arbitrarily small perturbations yield instantaneous amplification on $[0, T^\eps]$ to reach an $\mathcal{O}(1)$ amplitude.
 Note that $T^{\eps}$ converges to 0 as $\eps$ goes to 0.

Our result differs from the equivalent statement in \cite{Ge}: in (\ref{inst}), one obtains $\deltaZ\eps^{1/4}$ with the Dirichlet condition instead of $\deltaZ$ here, and therefore, the result with Dirichlet conditions is not a full instability result in the sense that it does not guarantee that perturbations reach an ${\cal O}(1)$ amplitude.
 This is due to the Prandtl boundary layers, the amplitude of which is the same as that of the internal layer. Indeed, after rescaling $(t,x,y)\mapsto\eps^{-1/2}(t,x,y)$ (these new coordinates are indifferently denoted $(t,x,y)$), a solution of $\eqref{NS(\sqrt{\eps})}$ with homogeneous Dirichlet boundary conditions can formally be written as (\ref{anzD}).
 In our case of a Navier condition depending on $\sqrt\eps$, the rescaled boundary condition no longer depends on the viscosity, so we will construct an approximate solution written as
$$ \ueps(t,x,y) \stackrel{\eps\rightarrow0}{\sim} \bust(y) + \eps^n \left(\uint(t,x,y) + \eps^{1/4}\ubl\left(t,x,\frac{y}{\eps^{1/4}}\right)\right) , $$
in which the factor $\eps^{1/4}$ in front of the boundary layer neutralises the one that appears when differentiating with respect to $y$. This ${\cal O}(\eps^{1/4})$ amplitude of the boundary layers is inspired by the asymptotic expansion shown by D.Iftimie and F.Sueur in \cite{IS}.
 In the Dirichlet case, this compensation does not take place, and considering times that are ${\cal O} \left(\left(n-\frac{1}{4}\right)\ln(\eps^{-1})\right)$ leads to a weaker form of instability.

The result is also different when $a=0$: it is shown in \cite{Bc} that, in this case, instability does not occur.
\newline

The proof of \thref{theo3} will follow Grenier's approach, which mainly relies on the construction of a WKB approximate solution to $(NS(a,\eps))$ starting with a linearly unstable solution to the Euler equation, and the use of a resolvent estimate on the linearised Euler equation (our \thref{resth}).
 The strength of the method comes from the fact that the resolvent estimate only needs to be shown on a subspace that does not have to be dense or of infinite dimension - typically spaces of functions with one or a finite number of Fourier modes.
 This has allowed it to be used in other circumstances than the study of boundary layers: F.Rousset and N.Tzvetkov used it to prove transverse nonlinear instability of solitary waves for KdV and water-wave equations (see \cite{RT1} and \cite{RT2}).

An example of smooth, linearly unstable boundary-layer profile, which can be considered as a model case, will be given after proof of \thref{theo3}.
\newline

\textbf{\underline{Organisation of the chapter.}} In the next section, we prove the $L^{2}$ convergence results. In section 3, we first prove the key resolvent estimate which will lead us to the proof of \thref{theo3}. In the final section, we provide an explicit profile that fit the instability result, from which we deduce a profile which will satisfy both theorems in the periodic setting.

\section{$L^{2}$ stability}

\subsection{Proof of \thref{theo1}}

We remind the reader of the notion of Leray solution to the Navier-Stokes system:

\begin{defi} $\ueps: \rplus\times\Omega\rightarrow \mathbb{R}^{2}$ is a Leray solution of $(NS(a,\eps))$ if:
\begin{enumerate}
\item $\ueps \in {\cal C}_{w}(\mathbb{R}^{+},L^{2}_{\sigma}) \cap L^{2}([0,T],H^{1}_{\sigma})$ for every $T>0$, where, if $E(\Omega)$ is a functional space on $\Omega$, $E_{\sigma}$ designates the space of divergence-free vector fields, tangent to the boundary and belonging to $E(\Omega)$.
\item $\ueps$ is a weak solution to $(NS(a,\eps))$ in the following sense: we have
$$ -\int_{\rplus}\int_{\Omega} \ueps\cdot \partial_{t}\varphi + 2a\sqrt{\eps}\int_{\rplus}\int_{\dOmega} \ueps\cdot\varphi + 2\eps\int_{\rplus}\int_{\Omega} S\ueps:S\varphi \hspace{50pt} $$
\begin{equation} \hspace{125pt}  - \int_{\rplus}\int_{\Omega} (\ueps\cdot\nabla\varphi)\cdot\ueps = \int_{\Omega} \ueps(0)\cdot\varphi(0) \label{WF} \end{equation}
for every $\varphi\in H^{1}(\rplus,H^{1}_{\sigma})$, where $A:B=\sum A_{i,j}B_{i,j}$ is the contracted product of two same-sized matrices $A$ and $B$.
\item For every $t\geq 0$, $\ueps$ satisfies the following energy estimate:
\begin{equation} \norme{L^{2}(\Omega)}{\ueps(t)}^{2}+4a\sqrt{\eps}\int_{0}^{t}\int_{\dOmega} |\ueps|^{2}+4\eps\int_{0}^{t}\int_{\Omega} |S\ueps|^{2} \leq \norme{L^{2}(\Omega)}{\ueps(0)}^{2} . \label{EE} \end{equation}
\end{enumerate}
\end{defi}

For $\eps>0$ fixed, such solutions are known to exist and are global in time. Moreover, in 2D, we have uniqueness. Likewise, the existence and uniqueness of a classical global-in-time solution $v$ to the Euler equation $\eqref{E}$ when $v_{0}\in H^{s}(\Omega)$, $s>2$, are well-known, and for every $T>0$, $v\in L^{\infty}([0,T],H^{s}(\Omega))$ (see \cite{Bc} or \cite{MB} and references therein).

Integration by parts leads to the following: let $f,~g\in H^{2}_{\sigma}$, with $f$ satisfying $(NC(a,\eps))$. Then
\begin{equation} -\int_{\Omega} \Delta f\cdot g = \frac{2a}{\sqrt{\eps}}\int_{\dOmega} f\cdot g + 2\int_{\Omega} Sf:Sg . \label{Lem} \end{equation}
This allows us to confirm that the weak formulation (\ref{WF}) contains the Navier-Stokes equation, the initial condition and the Navier boundary condition ($\eqref{NP}$ and the divergence-free conditions being given by the choice of the spaces in point 1 of the definition).
\newline

Let $\weps=\ueps-v$, which solves the equation
\begin{equation} \partial_{t}\weps + \ueps \cdot \nabla\weps+ \weps \cdot\nabla v- \eps\Delta\ueps +\nabla(\peps-q) = 0 . \label{6} \end{equation}
For $t>0$, $\weps$ satisfies the inequality
$$ \frac{1}{2}\norme{L^{2}(\Omega)}{\weps(t)}^{2} +\int_{0}^{t}\int_{\Omega}(\weps\cdot\nabla v)\cdot\weps + 2a\sqrt{\eps}\int_{0}^{t}\int_{\dOmega} \ueps\cdot\weps + 2\eps\int_{0}^{t}\int_{\Omega} S\ueps:S\weps \hspace{10pt} $$
\begin{equation} \hspace{200pt} \leq \frac{1}{2}\norme{L^{2}(\Omega)}{\weps(0)}^{2} , \label{7} \end{equation}
which can be seen as (\ref{6}) multiplied by $\weps$, integrated by parts in space and integrated in time.
 This direct computation is only possible because we are working in 2D and $\weps$ has the right regularity to be used as a test function, but (\ref{7}) also holds in 3D or higher; one starts by writing the energy inequality (\ref{EE}), then get three energy equalities: $(I)$ by multiplying the Euler equation by $v$, $(II)$ by again multiplying $\eqref{E}$ by $\ueps$, and $(III)$ by multiplying the Navier-Stokes equation by $v$.
 One then gets (\ref{7}) by doing $(\ref{EE})+(I)-(II)-(III)$ and using (\ref{Lem}). This idea appears, for example, in proofs of Kato-type theorems (\cite{Kt}, \cite{TW2}, \cite{Wx}).
\newline

For every $x,~y\in\mathbb{R}^{2}$, let $z=x-y$. We have
\begin{equation} 2(x\cdot z) = 2\left|z+\frac{y}{2}\right|^{2}-\frac{1}{2}|y|^{2} \label{AE} \end{equation}
and the same goes for the contracted matrix product. So, from (\ref{7}) we get
$$ \norme{L^{2}}{\weps(t)}^{2}+4a\sqrt{\eps}\int_{0}^{t}\int_{\dOmega} \left|\weps+\frac{v}{2}\right|^{2}+4\eps\int_{0}^{t}\int_{\Omega} \left|S\left(\weps+\frac{v}{2}\right)\right|^{2} \hspace{85pt} $$
\begin{equation} \hspace{65pt} \leq \norme{L^{2}}{\weps(0)}^{2}-2\int_{0}^{t}\int_{\Omega} (\weps\cdot\nabla v)\cdot\weps + a\sqrt{\eps}\int_{0}^{t}\int_{\dOmega} |v|^{2} + \eps\int_{0}^{t}\int_{\Omega} |Sv|^{2} . \label{13} \end{equation}
The left side is greater than $\norme{L^{2}}{\weps(t)}^{2}$, and we estimate each term on the right-hand side as follows:
\begin{itemize}
\item $\int_{\Omega}(\weps\cdot\nabla v)\cdot\weps \leq \norme{L^{\infty}}{\nabla v}\norme{L^{2}}{\weps}^{2}$ because $v(\tau)\in H^{s}(\Omega)$ for every $\tau\geq 0$ and for a certain $s>2$, so $v(\tau)\in{\cal C}^{1}(\bOmega)\cap H^{1}(\Omega)$, thus $\nabla v(\tau)\in L^{\infty}(\Omega)$;
\item $\int_{\dOmega}|v|^{2} \leq C_{\gamma}\norme{H^{1}}{v}^{2}$ for a certain $C_{\gamma}>0$;
\item $\norme{L^{2}}{Sv}^{2}\leq\norme{H^{1}}{v}^{2}$ given that $Sv=\frac{1}{2}(\nabla v+^{t}\nabla v)$;
\end{itemize}
so (\ref{13}) becomes
$$ \norme{L^{2}}{\weps(t)}^{2} \leq \norme{L^{2}}{\weps(0)}+\eps\left(\frac{aC_{\gamma}}{\sqrt{\eps}}+1\right)\int_{0}^{t}\norme{H^{1}}{v}^{2}~d\tau + 2\int_{0}^{t} \norme{L^{\infty}}{\nabla v}\norme{L^{2}}{\weps}^{2}~d\tau , $$
to which we apply Gr\"onwall's lemma, and we obtain
\begin{equation} \norme{L^{2}}{\weps(t)}^{2} \leq \left[\norme{L^{2}}{\weps(0)}+\eps\left(\frac{aC_{\gamma}}{\sqrt{\eps}}+1\right)\int_{0}^{t}\norme{H^{1}}{v}^{2}~d\tau\right]\exp\left(2\int_{0}^{t}\norme{L^{\infty}}{\nabla v}\right) . \label{14} \end{equation}

Now fix $T>0$, and take the supremum on the right-hand side for $t\in[0,T]$ to get a uniform-in-time bound on $\norme{L^2}{\weps(t)}$, thanks to which we conclude that
$$ \sup_{t\in[0,T]} \norme{L^{2}(\Omega)}{\weps(t)} = \mathop{{\cal O}}\limits_{\eps\rightarrow 0} \left(\max (\norme{L^2}{\weps(0)}, \eps^{1/4})\right). \hspace{30pt} \square $$

\subsection{Extension of the result (\thref{theostab})}

First, considering $a>0$, the entire proof of \thref{theo1} can be rewritten with $\eps^{-\beta}$ replacing $\eps^{-1/2}$ in the Navier condition, therefore replacing $\sqrt{\eps}$ by $\eps^{1-\beta}$ in the boundary terms of (\ref{13}). As a result, (\ref{14}) becomes
$$ \norme{L^{2}}{\weps(t)} \leq \left[\norme{L^{2}}{\weps(0)}+\eps^{1-\beta}C_v \left(aC_{\gamma}+ \eps^{\beta}\right)T\right]\exp\left(2C_v T\right) , $$
in which $C_v$ depends on norms of $v$, and the right-hand side clearly converges when $\beta<1$. However we cannot conclude when $\beta\geq 1$, because the coefficient $\eps^{1-\beta}aC_{\gamma}$ no longer goes to zero as $\eps\rightarrow0$.
\newline

We can also try rewriting the same proof with $a<0$ and the $(NC(a,\eps,\beta))$ boundary condition. Equality (\ref{7}) becomes
$$ \frac{1}{2}\norme{L^{2}(\Omega)}{\weps(t)}^{2} +\int_{0}^{t}\int_{\Omega}(\weps\cdot\nabla v)\cdot\weps + 2\eps\int_{0}^{t}\int_{\Omega} S\ueps:S\weps \hspace{70pt} $$
$$ \hspace{150pt} \leq \frac{1}{2}\norme{L^{2}(\Omega)}{\weps(0)}^{2} - 2a\eps^{1-\beta}\int_{0}^{t}\int_{\dOmega} \ueps\cdot\weps , $$
in which $-a>0$, thus the final boundary term cannot be ignored as in the positive case. Instead of equality (\ref{AE}), we use the more obvious $x\cdot z = |z|^2 + y\cdot z$ when $z=x-y$, so we get, after standard manipulations on the terms involving $v$:
$$ \frac{1}{2}\norme{L^{2}}{\weps(t)}^{2}+2\eps\int_{0}^{t}\norme{L^{2}}{S\weps}^{2} \leq \frac{1}{2}\norme{L^{2}}{\weps(0)}^{2}+\int_{0}^{t}\norme{L^{\infty}}{\nabla v}\norme{L^{2}}{\weps}^{2} $$
$$ \hspace{50pt} +(\eps+|a|\eps^{1-\beta})\int_{0}^{t}\norme{H^{1}}{v}^{2} + \eps\int_{0}^{t}\norme{L^{2}}{S\weps}^{2} + 2|a|\eps^{1-\beta}\int_{0}^{t}\int_{\dOmega}|\weps|^{2} . $$
We must deal with the troublesome final term $\Theta=2|a|\eps^{1-\beta}\norme{L^2(\dOmega)}{\weps}^2$. In what follows, $C$ is a constant that can change from one line to the next, but that never depends on $\eps$, $\beta$ or $\eta$, the latter being a parameter yet to appear.
\begin{itemize}
\item We start by using the trace theorem:
$$ \Theta\leq C\eps^{1-\beta}(\norme{L^{2}}{\weps}\norme{L^{2}}{\nabla \weps}+\norme{L^{2}}{\weps}^{2}); $$
\item we now use Young's inequality with a parameter $\eta$ to be chosen in a moment:
$$\Theta\leq C\eps^{1-\beta}\eta\norme{L^{2}}{\nabla \weps}^{2} + C\eps^{1-\beta}(\eta^{-1}+1)\norme{L^{2}}{\weps}^{2};$$
\item and finally, we use Korn's inequality, which in the half-space reads the following: for $f\in H^{1}_{\sigma}$, $\norme{L^{2}}{\nabla f}^{2}\leq 2\norme{L^{2}}{Sf}^{2}$. So, we get
$$ \Theta\leq C\eps^{1-\beta}\eta\norme{L^{2}}{S\weps}^{2} + C\eps^{1-\beta}(\eta^{-1}+1)\norme{L^{2}}{\weps}^{2} . $$
\end{itemize}
Thus we have made $\norme{L^{2}}{S\weps}^{2}$ appear - we therefore choose $\eta=\frac{1}{2C}\eps^{\beta}$ so that this term is absorbed by $\eps\norme{L^{2}}{S\weps}^{2}$ on the left-hand side. This leaves us with
$$ \norme{L^{2}}{\weps(t)}^{2} \leq \norme{L^{2}}{\weps(0)}^{2}+\int_{0}^{t} (\norme{L^{\infty}}{\nabla v} + C\eps^{1-\beta}(\eta^{-1}+1))\norme{L^{2}}{\weps}^{2} + C\eps^{1-\beta}(1+\eps^\beta)\int_{0}^{t}\norme{H^{1}}{v}^{2} , $$
which, after applying the Gr\"onwall lemma and choosing $t\leq T$ as in the proof above, leads to
$$ \norme{L^{2}}{\weps(t)}^{2} \leq \left[\norme{L^{2}}{\weps(0)}^{2}+C_{T}(\eps+\eps^{1-\beta})\right]e^{C_{T}(1+\eps^{1-\beta}+\eps^{1-\beta}\eta^{-1})} $$
for $C_{T}$ uniform in $t$, depending on norms of $v$, but not depending on $\eps$ or $\eta$. Noting that $\eta^{-1}\sim \eps^{-\beta}$, we have, for $\eps<1$,
$$ \norme{L^{2}}{\weps(t)}^{2} \leq \left[\norme{L^{2}}{\weps(0)}^{2}+C_{T}\eps^{1-\beta}\right]e^{C_{T}\eps^{1-2\beta}} , $$
which converges to zero if $\beta\leq 1/2$, because $e^{C\eps^{1-2\beta}}$ is unbounded for $\beta>1/2$. Notice the convergence rate in the first factor: ${\cal O}\left(\max(\norme{L^2}{\weps(0)}^2,\eps^{1-\beta})\right)$. $\square$

\section{Nonlinear instability}

As mentioned in the Introduction, a resolvent estimate on the Euler equation linearised around a shear flow $\ust$ will be required; we will begin by stating and proving it. To do so, we will consider individual Fourier modes and stream functions, and use the spectral theory of the linearised Euler equation in this setting with a fixed wave number, which is a one-dimensional problem.
 Then we shall prove \thref{theo3} by constructing an asymptotic expansion of solutions of $\eqref{NS(\eps)}$ around a linearly unstable shear profile $\ust$, as in \cite{Ge}. The terms in this WKB expansion will be wave packets, $U^{j}=\int_{\mathbb{R}} \varphi^{j}(k)V^{j}(k;t,y)e^{ikx}~dk$, where each $\varphi^{j}$ is smooth and compactly supported, thus $U^{j}$ is an $H^\infty$ function in the $x$-direction.
 We will choose $\varphi^{1}$ to be located around the most unstable wave number. Estimates on the important and unstable first term and the resolvent estimate will then give us the wanted instability. An example of linearly unstable shear flow will be given in the final paragraph.

\subsection{Preliminary results on the linearised Euler equation}

\subsubsection{Linear instability of the Euler equation}

Fix a wave number $k$, and consider the space $V_{k}$ of Fourier modes written as $u(t,x,y)=v(t,y)e^{ikx}$, with $v(t,\cdot)\in H^{s}$ for every $s>0$. We set $\sigma(k)$ the highest real part of complex numbers $\lambda$ that are unstable eigenvalues of the Euler equation linearised around $\ust$ on $V_k$; for each wave number $k$, it is finite, and $\sigma$ is an even analytic function of $k$.
 In turn, the function $\sigma$ has a maximum $\sigma_{0}>0$ (for details, see the study of the Rayleigh equation in \cite{Ge}, paragraph 4). We assume that this maximum is nondegenerate.

We write fields in $V_k$ that are tangent to the boundary and divergence-free by using stream functions:
\begin{equation} u(t,x,y)=e^{ikx}(\partial_{y}\Psi(t,y),-ik\Psi(t,y))=\nabla^{\bot}(e^{ikx}\Psi(t,y)) \label{stream} \end{equation}
with $\Psi(t,0)=0$. For $u\in V_{k}$, we use the norm $\norme{l}{u} = \sqrt{\norme{H^{l}_{k}}{u}^{2}+\norme{H^{l}_{k}}{\rot u}^{2}}$, where the $H^{l}_{k}$ norms are expressed as
$$ \norme{H^l_k}{u}^2 = \sum_{m=0}^l \norme{L^2}{\nabla^m_k v}^2 $$
with $\nabla ^m_k \Psi :=(k^{m_{1}}\partial^{m_{2}}_{y}\Psi(y))_{m_{1}+m_{2}=m}$.

From now on in the linear study, $u$ is in $V_k$ and written as in (\ref{stream}). Note that we have $\norme{H^l_k}{u}=\norme{H^l_k}{\nabla_{k}\Psi}$. As $\rot u = \partial_x u_2 - \partial_y u_1 = -\Delta(e^{ikx}\Psi(t,y)) = (k^{2}-\partial^{2}_{yy})\Psi$, standard elliptic regularity (see \cite{Kn} for example) gives us
\begin{equation} \norme{L^{2}_{k}}{\nabla u(t)} = \norme{L^{2}}{\nabla^{2}_{k}\Psi} \sim \norme{L^{2}_{k}}{\rot u(t)} . \label{450} \end{equation}
Using this, and standard properties of Sobolev spaces, we get that if $u\in V_{k}$ and $u'\in V_{k'}$, then $u\cdot\nabla u'\in V_{k+k'}$ and
\begin{equation} \norme{l}{u\cdot\nabla u'} \leq C\norme{l+2}{u}\norme{l+2}{u'} . \label{451} \end{equation}

\begin{theo} Let $\lambda'>\sigma_{0}$, and $w(t,x,y)\in V_{k}$ such that
\begin{equation} \norme{l}{w(t)}\leq C_{w}\frac{e^{\lambda't}}{(1+t)^{\alpha}} \label{reseq} \end{equation}
for every $l\geq 2$ and for some $\alpha\geq 0$. Consider the linearised Euler equation with source term
$$ \eqref{ELS}: \hspace{10pt} \left\{\begin{array}{rcl} \partial_{t}u+\ust\cdot\nabla u + u\cdot\nabla \ust + \nabla p & = & w \\
\div u & = & 0 \\ u|_{t=0} & = & u_{0} \in V_{k} \\ u_{2}|_{y=0} & = & u_{0,b}=e^{ikx}v_{0,b}(t) \end{array} \right. $$
with $u_{0,b}$ satisfying the estimate $|u_{0,b}(t)|\leq C_{w}(1+t)^{-\alpha}e^{\lambda't}$, and likewise for $\partial_{t}u_{0,b}$, and such that $u_{0,b}(0,x)=u_{0}(x)$. Then the solution of this system satisfies the estimate
$$ \norme{l-2}{u(t)}\leq C(1+t)^{-\alpha}e^{\lambda't} $$
for $t>0$, with $u(t)\in V_{k}$, and $C$, which depends on $w$, $\ust$, $u_{0}$, $u_{0,b}$, $l$ and $k$, and is locally bounded in the parameter $k$. \label{resth} \end{theo}

We show the result for functions with a single Fourier mode, but it extends to wave packets
$$ U(t,x,y) = \int_{\mathbb{R}} \varphi(k)u(k;t,x,y)~dk , $$
with $u(k)\in V_{k}$ and $\varphi$ a smooth compactly supported function such that, for every $k\in\supp(\varphi)$, $w(k)$ satisfies (\ref{reseq}) (see \cite{Ge} for example).
\newline

\preu{Proof:} we adapt the arguments used to prove the resolvent estimate for the KdV equation in \cite{RT1}. Set
$$ u(t,x,y) = e^{ikx}\vecdeux{\partial_{y}\Psi(t,y)}{-ik\Psi(t,y)} \hspace{10pt} \mathrm{and} \hspace{10pt} w(t,x,y) = e^{ikx}\tilde{w}(t,y) $$
and note that $\rot u = (k^{2}-\partial^{2}_{yy})\Psi:=B_{k}\Psi$, where, for $k\neq 0$, the differential operator $B_{k}:H^{l+2}(\rplus)\cap H^{1}_{0}(\rplus) \rightarrow H^{l}(\rplus)$ is invertible.

We first examine the $\norme{H^{l}_k}{\rot u}$ part of $\norme{l}{u}$. Taking the vorticity (rotational) of $\eqref{ELS}$, we obtain the equation
\begin{equation} \partial_{t}B_{k}\Psi + ik(\ust B_{k}\Psi+\ust''\Psi) = \rotk \tilde{w} , \label{40} \end{equation}
where $\rotk \tilde{w} = ik\tilde{w}_2-\partial_y\tilde{w}_1$. We then set
$$ \Psi_{0}(t,y) = \Psi(t,y) -e^{-\mu t}\Psi(0,y) -e^{-\mu y}\Psi(t,0) + e^{-\mu(t+y)}\Psi(0,0) $$
with an arbitrary $\mu>0$, thus $\Psi_{0}(0,y)=0$ and $\Psi_{0}(t,0)=0$. $\Psi_{0}$ solves
$$ \partial_{t}B_{k}\Psi_{0} + ik(\ust B_{k}\Psi_{0}+\ust''\Psi_{0}) = F_{0} . $$
We will not give the detailed expression of $F_{0}$, but we point out that integrating the $L^{2}$ hermitian dot-product of $F_{0}$ by $\Psi_{0}$ by parts leads to
$$ |(F_{0} | \Psi_{0})| \leq C\left(\norme{L^{2}}{f_{0}}^{2}+\norme{L^{2}}{\nabla_{k}\Psi_{0}}^{2} \right) $$
where $f_{0}$ contains terms depending on the data $w$, $u_{0}$ and $u_{0,b}$, and $\norme{L^{2}}{\nabla_{k}\Psi_{0}}$ is the norm of the velocity. More precisely, $f_{0} = \Phi_{1}+\Phi_{2}e^{-\mu t}$, where $\Phi_{1}$ contains $\rotk \tilde{w}$ and terms with $y=0$ (namely $u_{0,b}$ and $\partial_{t}u_{0,b}$), while $\Phi_{2}$ contains terms with $t=0$.

Using the Laplace transform in time, $\lapt g(z,y) = \int_{0}^{+\infty} e^{-zs}g(s,y)~ds$, and writing $\tilde{\Psi}=\lapt\Psi_{0}$ and $F=\lapt F_{0}$, we turn this differential equation into the eigenvalue problem:
\begin{equation} zB_{k}\tilde{\Psi} + ik(\ust B_{k}\tilde{\Psi}+\ust''\tilde\Psi) = F . \label{lapeq} \end{equation}
We choose $\gamma_{0}\in]\sigma_{0},\lambda'[$, and set $z=\gamma_{0}+i\tau$. $\gamma_{0}$ being fixed, we abbreviate $\tilde{\Psi}(z)=\tilde{\Psi}(\tau)$. When $\tau$ evolves in $\mathbb{R}$, we get the following estimates:

\begin{lemma} Let $l\geq 0$, and $\tilde{\Psi}$ solve (\ref{lapeq}) with $\tilde{\Psi}(\tau,0)=0$, and $F$ verifying
\begin{equation} |(F|\tilde{\Psi})| \leq C\left(\tilnorme{L^{2}}{\nabla_{k}\tilde{\Psi}}^{2} +\norme{L^{2}}{f}^{2} \right), \label{4050} \end{equation}
where $f=\lapt f_{0}$ depends on the data. There exists $C$ depending on $k$, $l$ and $\ust$, locally bounded in the parameter $k$, such that
$$ \tilnorme{H^{l}}{B_{k}\tilde{\Psi}(\tau)}^{2} \leq C\norme{H^{l+2}}{f(\tau)}^{2} \hspace{10pt} and \hspace{10pt} \tilnorme{H^{l}}{\nabla_{k}\tilde{\Psi}}^{2} \leq C\norme{H^{l+2}}{f}^{2} . $$ \label{laplem} \end{lemma}
Note that the lemma provides estimates on both the vorticity and the velocity. We prove this lemma in the next sub-paragraph.

By Parseval's equality and the above lemma, we have
\begin{eqnarray*} \int_{0}^{+\infty} e^{-2\gamma_{0}t}\norme{H^{l}}{B_{k}\Psi_{0}}^{2}~dt & = & \int_{-\infty}^{+\infty} \tilnorme{H^{l}}{B_{k}\tilde{\Psi}(\tau)}^{2}~d\tau \\
 & \leq & C\int_{-\infty}^{+\infty} \norme{H^{l+2}}{f(\tau)}^{2}~d\tau \\
 & \leq & C\int_{0}^{+\infty} e^{-2\gamma_{0}t}\norme{H^{l+2}}{\Phi_{1}(t)}^{2}~dt \\
 & & \hspace{25pt} + \int_{0}^{+\infty} e^{-2(\gamma_{0}+\mu)t}\norme{H^{l+2}}{\Phi_{2}}^{2} ~dt \\
 & \leq & C\int_{0}^{+\infty} e^{-2\gamma_{0}t}(\norme{H^{l+2}}{\Phi_{1}(t)}^{2}+\norme{H^{l+2}}{\Phi_{2}}^{2})~dt . \end{eqnarray*}
Replacing $\Phi_{1}$ by $\Phi_{1}\mathbf{1}_{[0,T]\times\rplus}$ for some $T>0$ does not affect the solution on $[0,T]\times\rplus$, so (\ref{reseq}) gives us
$$ \int_{0}^{T} e^{-2\gamma_{0}t}\norme{H^{l}}{B_{k}\Psi_{0}}^{2}~dt \leq C\int_{0}^{T} \frac{e^{2(\lambda'-\gamma_{0})t}}{(1+t)^{\alpha}}~dt , $$
where $C$ depends on $\norme{l}{u(0)}$, $\gamma_{0}$, $l$ and $k$, and is a locally bounded function of $k$. Noticing that
$$ \norme{H^{l}}{\Psi_{0}(t)} \geq \norme{H^{l}}{\Psi(t)} - \norme{H^{l}_{k}}{u_{0}}-\norme{H^{l}}{e^{-\mu y}}(|u_{0,b}(t)|+|\Psi(0,0)|) , $$
we get
\begin{equation} \int_{0}^{T} e^{-2\gamma_{0}t}\norme{H^{l}}{B_{k}\Psi}^{2}~dt \leq C\int_{0}^{T} \frac{e^{2(\lambda'-\gamma_{0})t}}{(1+t)^{\alpha}}~dt \leq C\frac{e^{2(\lambda'-\gamma_{0})T}}{(1+T)^{\alpha}} \label{403} \end{equation}
(the last inequality is obtained by integrating by parts).

Using the same procedure, estimate (\ref{403}) also holds for the $\norme{H^{l}_k}{u}$ part of $\norme{l}{u}$:
\begin{equation} \int_{0}^{T} e^{-2\gamma_{0}t}(\norme{H^{l}}{k\Psi}^{2} + \norme{H^{l}}{\partial_{y}\Psi}^{2})~dt \leq C\frac{e^{2(\lambda'-\gamma_{0})T}}{(1+T)^{\alpha}} . \label{404} \end{equation}

A quick energy estimate on (\ref{40}) gives us
$$ \frac{1}{2}\partial_{t}\norme{H^{l}}{B_{k}\Psi(t)}^{2} \leq C(\norme{H^{l}}{B_{k}\Psi(t)}^{2}+\norme{H^{l}_{k}}{u(t)}^{2}+\norme{H^{l}_{k}}{\rot w(t)}^{2}) . $$
Multiply this by $e^{-2\gamma_{0}t}$, and using the hypothesis on $w$, we get
$$ \partial_{t}(e^{-2\gamma_{0}t}\norme{H^{l}}{B_{k}\Psi}^{2}) \leq C\left(e^{-2\gamma_{0}t}\norme{H^{l}}{B_{k}\Psi}^{2}+e^{-2\gamma_{0}t}\norme{H^{l}_{k}}{u}^{2}+\frac{e^{2(\lambda'-\gamma_{0})t}}{(1+t)^{\alpha}}\right) . $$
Finally, we integrate between 0 and $T$, and use (\ref{403}) and (\ref{404}) to we get the desired estimate: $\norme{H^{l}_{k}}{\rot u(t)} \leq C(1+t)^{-\alpha}e^{\lambda't}$.

Thanks to (\ref{450}), it only remains to estimate $\norme{L^{2}_k}{u}$. We simply multiply the linearised Euler equation $\eqref{ELS}$ by $u$, and get
$$ \frac{1}{2}\partial_{t}\norme{L^{2}_{k}}{u}^{2} \leq C(\norme{L^{2}_{k}}{u}^{2} + \norme{L^{2}_{k}}{\nabla u}^{2} + \norme{L^{2}_{k}}{w}^{2}), $$
where $C$ depends only on $\ust$. Multiplying by $e^{-2\gamma_{0}t}$, integrating in time and combining (\ref{450}), (\ref{403}) and (\ref{404}), we get $\norme{L^{2}_{k}}{u(t)}\leq C(1+t)^{-\alpha}e^{\lambda't}$. $\square$

\subsubsection{Proof of the resolvent lemma (\lemref{laplem})}

First note that the estimate is obvious when $k=0$. We must then prove the estimate for $k\neq 0$ and then show that the constant $C$ is a locally bounded function of $k$ as $k\rightarrow 0$.
 Therefore, we consider $|k|\leq K$, and divide the proof in two parts: with $z=\gamma_{0}+i\tau$ and $\tau\in \mathbb{R}$, we show the estimate for $|\tau|\geq M$, with $M$ to be chosen independent of $l$, and then for $|\tau|\leq M$, using the notion of exponential dichotomy (see \cite{Cw}).
\newline

\textbf{Step 1: $|\tau|$ large.} We start with an $L^{2}$ estimate. Set $\Theta=B_k\tilde{\Psi}$. The imaginary part of the dot-product of (\ref{lapeq}) with $\Theta$ leads to
\begin{equation} |\tau|\norme{L^{2}}{\Theta}^{2} \leq (C_{sh}|k|+1)\norme{L^{2}}{\Theta}^{2} + \frac{1}{2}\norme{L^{2}}{F}^{2} + C_{sh}\tilnorme{L^{2}}{\nabla_{k}\tilde{\Psi}}^{2} , \label{550} \end{equation}
noting that $\tilnorme{L^{2}}{k\tilde{\Psi}}^{2}\leq \tilnorme{L^{2}}{\nabla_{k}\tilde{\Psi}}^{2}$. To deal with this final term, we consider the real part of the dot-product of (\ref{lapeq}) and $\tilde{\Psi}$: we integrate $(\partial_{yy}^2 \tilde{\Psi}|\tilde{\Psi})$ by parts, notice that $(u_s '' \tilde{\Psi}|\tilde{\Psi})$ is real and use inequality (\ref{4050}) as well as Young's inequality with parameter $\gamma_0$, to get
$$ \gamma_{0}\tilnorme{L^{2}}{\nabla_{k}\tilde{\Psi}}^{2} \leq \left[C_{sh}|k|+\frac{\gamma_{0}}{2}\right]\tilnorme{L^{2}}{\nabla_{k}\tilde{\Psi}}^{2}+C\norme{H^{1}}{f}^{2} , $$
having noticed that $\norme{L^2}{F}^2\leq\norme{H^1}{f}^2$. Let $k_0>0$ be such that $C_{sh} k_{0}+\frac{\gamma_{0}}{2}=\frac{3\gamma_{0}}{4}$. When $|k|\leq k_0$, we can absorb the first term on the right-hand side to get $\tilnorme{L^{2}}{\nabla_{k}\tilde{\Psi}}^{2}\leq C\norme{H^{1}}{f}^{2}$, with $C$ depending only on $\gamma_{0}$ and $\ust$. This allows us to conclude from estimate (\ref{550}):
$$ h(\tau,k)\norme{L^{2}}{\Theta}^{2} = (|\tau|-C_{sh}|k|-1)\norme{L^{2}}{\Theta}^{2} \leq C\norme{H^{1}}{f}^{2}. $$
By choosing $M>0$ such that $h(M,k_0)\geq 1$, we have the desired uniform estimate for $|\tau|\geq M$ and $|k|\leq k_0$. In particular, we have the desired local boundedness of the constant $C$ when $|k|\rightarrow 0$.

When $K\geq |k|\geq k_{0}$, things are easier, as we have uniform elliptic regularity for the operator $B_k$. Consider the imaginary part of the $L^{2}$ dot-product of equation (\ref{lapeq}) with $\Theta=B_{k}\tilde{\Psi}$: thanks to the elliptic regularity bound $\norme{H^{2}}{\Psi}\leq \max(1,k_{0}^{-2})\norme{L^{2}}{B_{k}\Psi}$, we have
$$ |\tau|\norme{L^{2}}{\Theta}^{2} \leq \left[C_{sh}|k|+\frac{1}{2}\right]\norme{L^{2}}{\Theta}^{2}+\frac{1}{2}\norme{L^{2}}{F}^{2} , $$
so we can write once again $h(\tau,k)\norme{L^2}{\Theta}^2 \leq C\norme{H^1}{f}^2$, with $C$ depending on $\ust$ and $k_{0}$. Choosing $M$ so that $h(M,K)\geq 1$ (remember that $|k|$ is assumed to be bounded), we have, for $|\tau|\geq M$, $\norme{L^{2}}{\Theta}\leq C\norme{H^{1}}{f} \leq C\norme{L^{2}}{f}$. This ends the proof of the $L^2$ estimates for $|\tau|$ large.
\newline

$H^{l}$ estimates are easily obtained by induction: examine the dot-product of $\partial^{l}$(\ref{lapeq}) by $\partial^{l}\Theta$, and use Young's inequality and the induction hypothesis $\norme{H^{l-1}}{\Theta}\leq C\norme{H^{l-1}}{F}$ to get $\frac{\gamma_{0}}{2}\norme{L^{2}}{\partial^{l}\Theta}^{2}\leq C\norme{H^{l}}{F}^{2}$, with $C$ depending on $\gamma_{0}$, $|k|$ and $l$.
\newline

\textbf{Step 2: $|\tau|$ small.} We seek a solution to (\ref{lapeq}) written as $\tilde{\Psi}=\Psi_{1}+\Psi_{2}$, where
\begin{eqnarray} zB_{k}\Psi_{1} & = & F , \label{410} \\ \mathrm{and} \hspace{10pt} (z+ik\ust)B_{k}\Psi_{2}+ik\ust''\Psi_{2} & = & -ik(\ust B_{k}\Psi_{1}+\ust''\Psi_{1}) := kG \label{411} \end{eqnarray}
Immediately, we get $\norme{H^{l}}{B_{k}\Psi_{1}} \leq \frac{1}{\gamma_{0}}\norme{H^{l}}{F}$, and, by multiplying (\ref{410}) by $\Psi_{1}$, we get
$$ \frac{\gamma_{0}}{2}\norme{L^{2}}{\nabla_{k}\Psi_{1}}^{2} \leq C\norme{L^{2}}{f}^{2} , $$
(it is proved similarly to (\ref{4050}), but also using Young's inequality) and likewise with $H^l$ norms, so a $\norme{H^{l}}{\Psi_{2}}\leq C\norme{H^{l}}{kG}$ estimate suffices to prove the result, as
$$ \norme{H^{l}}{kG} \leq c|k|\norme{H^{l}}{B_{k}\Psi_{1}}+c'\norme{H^{l}}{\nabla_{k}\Psi_{1}} . $$
First rewrite (\ref{411}) as an ordinary differential system:
\begin{equation} \partial_{y} U(y) = A(k,z;y)U(y) + H(y) \label{412} \end{equation}
with $U = (\Psi_{2},\partial_{y}\Psi_{2})$, $H=(0,-(z+iku_s)^{-1}kG)$ (we have $|H|\leq \gamma_0^{-1} |kG|$), and
$$ A(k,z;y) = \matdeux{0}{1}{k^{2}+\frac{ik\ust''(y)}{z+ik\ust(y)}}{0} = A_{\infty}(k) + B(k,z;y) , $$
where $A_{\infty}(k)=\lim_{y\rightarrow+\infty}A(k,z;y)$, and $B(k,z;y)={\cal O}(e^{-\eta y})$ for some $\eta>0$ since $\ust''$ decays exponentially.
 The eigenvalues of $A_{\infty}(k)$ are real, so by the roughness of exponential dichotomy \cite{Cw}, the system $\partial_{y}U=A(k,z)U$ has an exponential dichotomy on $\rplus$; this means that if $T(k,z;y,y')$ is the fundamental solution of this last equation with $T(k,z;y,y)=I_{2}$, there exist a projection $P(k,z;y)$, verifying
$$ T(k,z;y,y')P(k,z;y')=P(k,z;y)T(k,z;y,y') , $$
and positive constants $C$ and $\alpha$, with all of these depending smoothly on $(k,z)$, such that for any $\xi\in\mathbb{C}^{2}$,
$$ \begin{array}{rcll} |T(k,z;y,y')P(k,z;y')\xi| & \leq & C(k,z)e^{-\alpha(y-y')}|\xi| & \mathrm{if}~y\geq y'\geq 0 \\
|T(k,z;y,y')(I-P(k,z;y'))\xi| & \leq & C(k,z)e^{\alpha(y-y')}|\xi| & \mathrm{if}~y'\geq y\geq 0 . \end{array} $$
For any $0<\rho<\rho'$, $C(k,z)$ is uniformly bounded in the set
$$ {\cal K}_{\rho,\rho'} = \{(k,\gamma_{0}+i\tau)~|~|k|\in [\rho,\rho']~\mathrm{and}~|\tau|\leq M\},  $$
but we would like $C(k,z)$ to be bounded in ${\cal K}_{0,\rho'}$ for a certain $\rho'$.
\newline

To get that $C(k,z)$ is bounded near $k=0$, we follow the proof of the persistence of ordinary dichotomies as done in \cite{Cw}. $A_{\infty}(k)$ cannot be uniformly diagonalised near $k=0$, because the basis of diagonalisation for $k\neq 0$ is $\matdeux{1}{1}{-|k|}{|k|}$.
 Instead, we change to the basis $(v_{1}(k), v_{2}(k))=\matdeux{1}{0}{-|k|}{1}$. $v_{1}$ is an eigenvector of $A_{\infty}(k)$, with eigenvalue $-|k|$. $A_{\infty}(k)$ is therefore triangular in this basis, and $v_{1}$ spans a space corresponding to that of exponentially decreasing solutions of the equation $\partial U=A_{\infty}(k)U$.
 Setting $T_\infty$ the fundamental solution of that equation, there exists a projection $\Pi$ such that
\begin{equation} \begin{array}{rcll} |T_{\infty}(k;y,y')(I-\Pi(k;y'))\xi| & \leq & e^{-|k|(y-y')}|\xi| & \mathrm{if}~y\geq y'\geq 0 \\
|T_{\infty}(k;y,y')\Pi(k;y')\xi| & \leq & e^{|k|(y-y')}|\xi| & \mathrm{if}~y'\geq y\geq 0 . \end{array} \label{413} \end{equation}
Consider $E(k)$, the Banach space of functions $V$ such that
$$ \norme{E(k)}{V}:=\tilnorme{{\cal C}^{0}}{Ve^{|k|y}}<+\infty , $$
and $S$ the linear mapping defined by
$$ SV(y) = \int_{0}^{y} T_{\infty}(k;y,y')(I-\Pi(k;y'))B(k,z;y')V(y')~dy' \hspace{45pt} $$
$$ \hspace{45pt} -\int_{y}^{+\infty} T_{\infty}(k;y,y')\Pi(k;y')B(k,z;y')V(y')~dy' . $$
Remember that $|B(k,z;y)|\leq b|k|e^{-\eta y}$. Using (\ref{413}), we get that if $V\in E(k)$, then $SV\in E(k)$ with the estimate
\begin{equation} |SV(y)e^{|k|y}| \leq \norme{E(k)}{V}\left(\int_{0}^{+\infty} b|k|e^{-\eta y'}~dy'\right) \leq \frac{1}{2}\norme{E(k)}{V} \label{415} \end{equation}
when $|k|\leq \rho'$ small enough. So $S$ is a contracting endomorphism of $E(k)$.

Let $U_{\infty}(k)\in E(k)$ be a solution of $\partial U=A_{\infty}(k)U$. By the Duhamel formula, a bounded solution of the equation $\partial U = AU$ is a fixed point of the affine transform $\tilde{S}=S+U_{\infty}$; thanks to (\ref{415}), Picard's fixed point theorem allows us to conclude that such a solution $U$ exists in $E(k)$. Finally, we must uniformly bound $\norme{E(k)}{U}$. Choosing $U_{\infty}(k)=(e^{-|k|y},-|k|e^{-|k|y})$, the decreasing eigenfunction of $A_{\infty}$, (\ref{415}) gives us the wanted bound.
\newline

The end of the proof is the same as in \cite{RT1}; we provide it for completeness. Note that for $k\neq 0$, $I-P(k,z): f \mapsto [y\mapsto (I-P(k,z;y))f(y)]$ is the projection on the subspace of solutions that go to 0 when $y\rightarrow+\infty$; let us define $Q(k,z)$ the projection on the subspace of solutions of which the first component vanishes at $y=0$.
 As the linearised Euler equation does not have an eigenvalue with real part $\gamma_{0}$, necessarily we have
\begin{equation} {\cal R}(I-P(k,z;0))\cap{\cal R}(Q(k,z;0)) = \{0\} , \label{supp} \end{equation}
where we denote ${\cal R}(A)$ the range of a matrix $A$, and we define a basis of solutions $(e_{1}(k,z),e_{2}(k,z))$, with $e_{1}\in{\cal R}(I-P(k,z;0))$ and $e_{2}\in{\cal R}(Q(k,z;0))$. Define new projections
$$ P'(k,z) = (e_{1}(k,z),e_{2}(k,z)) \matdeux{1}{0}{0}{0} (e_{1}(k,z),e_{2}(k,z))^{-1} , $$
and $P'(k,z;y)=T(k,z;y,0)P'(k,z)$, so that we have both ${\cal R}(P'(k,z;y))={\cal R}(Q(k,z;y))$ and ${\cal R}(I-P'(k,z;y))={\cal R}(I-P(k,z;y))$. We also have the estimates
\begin{equation} \begin{array}{rcll} |T(k,z;y,y')P'(k,z;y')\xi| & \leq & C'(k,z)e^{-\alpha(y-y')}|\xi| & \mathrm{if}~y\geq y'\geq 0 \\
|T(k,z;y,y')(I-P'(k,z;y'))\xi| & \leq & C'(k,z)e^{\alpha(y-y')}|\xi| & \mathrm{if}~y'\geq y\geq 0 . \end{array} \label{pest} \end{equation}
Again, we must check that $C'(k,z)$ is bounded in ${\cal K}_{0,\rho}$. To do so, we point out that the projections $Q(k,z)$ and $I-P(k,z)$ can be continued up to $k=0$: we have that ${\cal R}(Q(0,z))$ is the subspace of solutions to the equation $z\Psi''=0$ with $\Psi(0)=0$, and ${\cal R}(I-P(0,z))$ is the subspace of bounded solutions to $z\Psi''=0$ (constants).
 As the only bounded solution of $z\Psi''=0$ with $\Psi(0)=0$ is $\Psi\equiv0$, (\ref{supp}) is true up to $k=0$, so $C'(k,z)$ is bounded for $k$ near 0.

Finally, by Duhamel's formula, a bounded solution of (\ref{412}) is, for fixed $(k,z)$,
$$ U(y) = \int_{0}^{y} T(y,y')P'(y')H(y')~dy' - \int_{y}^{+\infty} T(y,y')(I-P(y'))H(y')~dy' , $$
as the only solution of $\partial_y U = A(k,z)U$ with $U_1\in H^1_0(\rplus)$ is zero. So, by using (\ref{pest}) and standard convolution estimates, we have $\norme{L^{2}}{U} \leq C(k,z)\norme{L^{2}}{H}$. To get estimates on the $L^{2}$ norms of the derivatives, just differentiate the equation $l$ times, notice that it is the same type of system and conclude by induction. $\square$

\subsection{Proof of the instability result}

Recall that $\ueps$ is a solution of $\eqref{NS(a,\eps)}$. First and foremost, as the initial profile depends on $\eps^{-1/2}y$, we rescale the variables $(t,x,y)\rightarrow \eps^{1/2}(t,x,y)$, and these new variables will be indifferently noted $(t,x,y)$. From now on, the faster variable will be $Y=\eps^{-1/4}y$. The system we get after rescaling is
$$ \eqref{NS'}: \hspace{10pt} \left\{ \begin{array}{rrcl} \eqref{NS(\sqrt{\eps})}: \hspace{10pt} & \partial_{t}\ueps+\ueps\cdot\nabla\ueps-\sqrt{\eps}\Delta\ueps+\nabla\peps & = & 0 \\
& \div\ueps & = & 0 \\
& \ueps(0,x,y) & = & \ueps_{0}(x,y) \\
\eqref{NP}: \hspace{10pt} & \ueps_{2}(t,x,0) & = & 0 \\
(NC(a)): \hspace{10pt} & \frac{1}{2}\partial_{y}\ueps_{1}(t,x,0) & = & a\ueps_{1}(t,x,0) \end{array} \right. $$
with the initial condition $\ueps_{0}$ close to a linearly unstable shear flow $\ust$. The Navier boundary condition no longer depends on the viscosity, and we will show the following asymptotic expansion:
\begin{equation} \ueps(t,x,y) \sim \bust(t,y) + u^{i}(t,x,y) + \eps^{1/4}u^{b}(t,x,Y) , \label{ISAE} \end{equation}
where $u^{i}$ solves the Euler equation and $u^{b}$ is the boundary layer.

Starting with a linearly unstable solution to the Euler equation, we will construct an approximate solution to the Navier-Stokes system $\eqref{NS'}$ above that will allow us to prove the instability inequality (\ref{inst}). The approximate solution is described in the following proposition.

\begin{propo} For given integers $n$ and $N$, there exist a time $T^\eps_0$, depending on $\eps$ but not converging to zero, a constant $C_R>0$ and a family $((U^j,P^j))_{j\in \{1,\cdots,N\}}$ of functions such that the vector field
$$ \uap(t,x,y) = \bust(t,y) + \sum_{j=1}^N \eps^{jn}U^j(t,x,y) $$
is divergence-free and tangent to the boundary, and the pair $(\uap,\pap)$, with $\pap$ also written as a WKB expansion, $\pap=\sum_{j=1}^N \eps^{jn}P^j$, solves the Navier-Stokes equation approximately, in the sense that $\norme{L^2(\Omega)}{\uap(0)-\ust} \leq C\eps^n$, and
$$\partial_t \uap + \uap \cdot\nabla \uap - \sqrt\eps \Delta\uap + \nabla\pap = R, $$
with $R$ satisfying the following growth bound on the time interval $[0,T^\eps_0]$:
\begin{equation} \norme{L^2}{R(t)} \leq \frac{C_R\eps^{(N+1)n}}{(1+t)^{(N+1)/2}}(1+t)^{1/4}e^{(N+1)\sigma_{0}t}. \label{Rest} \end{equation}

Moreover, there exist constants $C_j$, $j\in\{1,\cdots,N\}$ such that the individual components of $\uap$ satisfy the following inequalities on $[0,T^\eps_0]$:
\begin{eqnarray} \norme{H^{1}}{U^{j}(t)} & \leq & C_{j}(1+t)^{1/4}\frac{e^{j\sigma_{0}t}}{(1+t)^{j/2}}, \label{est1} \\
\norme{L^\infty}{U^{j}(t)} & \leq & \frac{C_je^{j\sigma_{0}t}}{(1+t)^{j/2}}, \label{estinf} \end{eqnarray}
and there exists $C'_1>0$ and a bounded domain $\Omega_A (t)\subset\Omega$, whose measure is of order $\sqrt{1+t}$, on which
\begin{equation} \norme{L^{2}(\Omega_{A}(t))}{U^{1}(t)} \geq C'_{1}(1+t)^{1/4}\frac{e^{\sigma_{0}t}}{\sqrt{1+t}}. \label{est2} \end{equation}
\label{approx} \end{propo}

This proposition will be proved over the next two sections: the first to explain the construction of $\uap$, in which the role of $u^\sharp$ will be given in particular, and the second to prove the estimates. The final estimate in \propref{approx} will allow us to get a lower bound on the $L^\infty$ norm of $U^1$, which will be crucial to prove the instability in the final paragraph.

\subsubsection{Proof of \propref{approx}, part 1 \\ Building an approximate solution}

For now, fix $N$ large and arbitrary (it will be chosen in section 3.2.3). Each $U^{j}$ will be written as a wave packet; letting $k_{0}$ be such that $\sigma(k_{0})=\sigma_{0}$, and setting $\varphi^{j}$ as compactly supported smooth functions of $k$ with $k_{0}\in\mathrm{supp}(\varphi)$, we write
$$ U^{j}(t,x,y) = \int_{\mathbb{R}}\varphi^{j}(k)V^{j}(k;t,x,y)~dk , $$
with the $V^{j}(k)$ being in the previously-defined $V_{k}$ space. We shall write these wave packets more precisely in the next paragraph.

The equations that the $U^{j}$ are supposed to solve are Navier-Stokes equations linearised around $\bust$:
\begin{equation} \partial_{t}U^{j} + \bust\cdot\nabla U^{j} + U^{j}\cdot\nabla\bust - \sqrt\eps \Delta U^{j} +\nabla P^{j} + \sum_{j_{1}+j_{2}=j} U^{j_{1}}\cdot\nabla U^{j_{2}} = 0 . \label{nsl} \end{equation}
We will not necessarily have $\div U^{j}=0$ individually, but in total $\div \uap$ must be zero. We again only solve (\ref{nsl}) approximately, by writing a sub-expansion
\begin{equation} U^{j} = \sum_{m=0}^{8n-1} \eps^{m/8} u^{i,8(j-1)n+m}(t,x,y)+\eps^{(m+2)/8} u^{b,8(j-1)n+m}(t,x,Y) , \label{anzl2} \end{equation}
and likewise for the pressure - there is a gap of $\eps^{1/4}$ between corresponding internal and boundary layers, which is inspired by the asymptotic expansion (\ref{ISAE}). We therefore solve (\ref{nsl}) with an error $\eps^{n}{\cal E}^j$, which is taken into account in the equation on $U^{j+1}$, so we solve, at each level,
$$ \partial_{t}U^{j} + \bust\cdot\nabla U^{j} + U^{j}\cdot\nabla\bust -\sqrt\eps \Delta U^{j} + \nabla P^j = -{\cal E}^{j-1} - \sum_{j_{1}+j_{2}=j} U^{j_{1}}\cdot\nabla U^{j_{2}} . $$
In the final term $U^N$, we will need to construct more than $8n$ terms, to ensure that the terms in the error ${\cal E}^N$ will have the right growth. This also means that the final error can be as small as we want, of order $\eps^{n+l_0}$, by choosing $l_0$ large. $U^N$ will also contain a corrective term $u^\sharp$ to ensure that $\div\uap = 0$ and $\uap_2|_{y=0}=0$.
\newline

Writing $q=8(j-1)n+m$, we must now understand the equations that the $u^{i,q}$ and $u^{b,q}$ solve to get estimates on these functions.

We plug ansatz (\ref{anzl2}) into each equation of the linearised Navier-Stokes system, and consider the equations obtained when asking terms of a same order $\eps^{n+q/8}$ to cancel out. We get that the internal layers $u^{i,q}$, which solve only the Euler part of the $\eqref{NS(\sqrt{\eps})}$ equation, are solutions of linear Euler systems
$$ (EL(q)): \hspace{10pt} \left\{ \begin{array}{rrl} (E_{q}): \hspace{10pt} & \partial_{t}u^{i,q} + \ust\cdot\nabla u^{i,q}+u^{i,q}\cdot\nabla\ust + r^{i,q} + \nabla p^{i,q} & = 0 \\ & \div u^{i,q} & = 0 \\ (P_{q}): \hspace{10pt} & u^{i,q}_{2}(t,x,0) + u^{b,q-2}_{2}(t,x,0) & =0 \end{array} \right. $$
where we linearise around $\ust$ instead of $\bust$, because we have a better understanding of the linear equation around time-independent profiles, and
$$ r^{i,q} = \frac{\bust-\ust}{\eps^{1/8}}\cdot\nabla u^{i,q-1} + u^{i,q-1}\cdot\nabla\left(\frac{\bust-\ust}{\eps^{1/8}}\right) - \Delta u^{i,q-4} + \sum_{8n+q_{1}+q_{2}=q} u^{i,q_{1}}\cdot\nabla u^{i,q_{2}} . $$
We also get that the boundary layers each solve a Stokes problem with a Neumann boundary condition
$$ (S(q)): \hspace{10pt} \left\{ \begin{array}{rrl} (S_{q}): \hspace{10pt} & \partial_{t}u^{b,q} - \partial^{2}_{YY}u^{b,q} + r^{b,q}+\nabla p^{b,q} & = 0 \\
(D_{q}): \hspace{10pt} & \partial_{Y}u_{2}^{b,q} + \partial_{x}u_{1}^{b,q-2} & = 0 \\
(N_{q}): \hspace{10pt} & \left(\frac{1}{2} \partial_{Y}u_{1}^{b,q} + \frac{1}{2}\partial_{y}u^{i,q}_{1}-au^{i,q}_{1}-au_{1}^{b,q-2}\right)(t,x,0) & = 0 \\
(V_{q}): \hspace{10pt} & \lim_{Y\rightarrow +\infty} u_{2}^{b,q}(t,x,Y) & = 0.
\end{array} \right. $$
We do not give any detail on $r^{b,q}$ other than it depends on $\eps^{-1/8}\bust$ and terms $u^{i,m}$ and $u^{b,m}$ with $m<q$.

In both systems, initial conditions are chosen to be compatible with the boundary conditions, and rapidly decreasing functions of $y$ (like $v(0,x,0)e^{-y}$); therefore $U^{j}|_{t=0}$ is ${\cal O}(\eps^{n})$. For low values of $q$, any term with a negative index is of course ignored.
\newline

Now that we know the equations that each layer is supposed to solve, we can construct the approximate solution by induction. As announced, we start by choosing $u^{i,0}$ as an unstable solution of the linear Euler equation $(E_{1})$ with homogeneous non-penetration condition. We write
$$ u^{i,0}(t,x,y) = \int_{\mathbb{R}} \varphi^{1}(k)v^{i,0}(k;t,x,y)~dk , $$
for a compactly supported function $\varphi^1$ to be chosen a little later, and we choose the most unstable mode for each $k$, \textit{i.e.}
$$ v^{i,0}(k;t,x,y) = e^{\lambda(k)t}\tilde{v}(k;y)e^{ikx} , $$
with $\Re(\lambda(k))=\sigma(k)$. Then, to solve $(S(0))$, $(D_{0})$ and the limit condition allow us to choose $u^{b,0}_{2}=0$, then $p^{b,0}=0$ and $u^{b,0}_{1}$ solves a heat equation.

Having built $u^{i,m}$ and $u^{b,m}$ for $m<q$, $u^{i,q}$ is the solution of the linear Euler system $(E(q))$. Now we build $u^{b,q}$ by taking the equations of the Stokes system $(S(q))$ one after another.
 First, $(D_{q})$ is the equation relative to the divergence-free condition: integrating it between $Y$ and $+\infty$ (to satisfy $(V_q)$), we get an expression of the normal component of $u^{b,q}$ that goes to 0 as $Y\rightarrow+\infty$; the rapid-decrease property will follow from the induction below, and nothing guarantees that we have $u^{b,q}_{2}(t,x,0)=0$, so this is where the boundary condition $(P_{q})$ in the Euler system comes from.
 Then, notice that thanks to the structure of $\bust$, the normal part of the Stokes equation $(S_{q})$ depends only on $u^{b,q}_{2}$ and $p^{b,q}$ (and of course previous terms), so we get $p^{b,q}$ by the same method.
 Finally, the first component of $(S_{q})$ is now a heat equation, with the Neumann boundary condition $(N_{q})$. We continue the construction up to $M=8nN+8l_0$; all terms with $q\geq 8nN$ are part of $U^N$.
\newline

We now consider a corrective term $\eps^{n+l_0}u^\sharp(t,x,y)$ (which is part of $U^N$ as mentioned above) such that
\begin{eqnarray} \div u^\sharp & = & -\partial_{x}v^\sharp _1 \label{coreq} \\
\mathrm{and} \hspace{10pt} u^\sharp_{2}(t,x,0) & = & -v^\sharp _2(t,x,0), \label{corbd} \end{eqnarray}
where $v^\sharp = u^{b,M-1}+\eps^{1/8}u^{b,M}$. This ensures that $\uap$ is divergence-free and tangent to the boundary, although each $U^j$, taken individually, is not. Note that, because of $(D_q)$, $\int_{\Omega} \partial_x v^\sharp_1 = \int_{\dOmega} v^\sharp_2 = 0$, so the equation is consistent with the divergence theorem.

To construct a solution, we decompose
$$ u^\sharp(x,y)=u^{\sharp,1}(x,y)+u^{\sharp,2}(x,\eps^{-1/4}y), $$
so that $u^{\sharp,1}$ is divergence-free with $u^{\sharp,1}_2(t,x,0)=-v^\sharp_2(t,x,0)$. We simply set
$$ u^{\sharp,1}_2(t,x,y)=-v^\sharp_2(t,x,0)\chi(y) $$
with $\chi$ a smooth function supported in $[0,1]$ with $\chi(0)=1$, and get $u^{\sharp,1}_1$ by integrating the divergence-free condition and using $(D_M)$ and $(D_{M-1})$:
\begin{eqnarray*} u^{\sharp,1}_1(t,x,y) & = & \chi '(y) \int_{x}^{+\infty} v_2^\sharp(t,\xi,0)~d\xi \\
 & = & - \chi '(y)\int_x^{+\infty} \int_0^{+\infty} \partial_x (u^{b,M-3}_1+\eps^{1/8}u^{b,M-2}_1)(t,\xi,Y)~d\xi~dY \\
 & = & \chi '(y) \int_{0}^{+\infty} u^{b,M-3}_1(t,x,Y) + \eps^{1/8}u^{b,M-2}_1(t,x,Y)~dY. \end{eqnarray*}
Meanwhile, $u^{\sharp,2}$ will be tangent to the boundary and $\div u^{\sharp,2} = -\partial_x v^\sharp_1$. This equation can be solved classically (as in the introduction of \cite{ASV} for example), by decomposing $u^{\sharp,2}=\nabla \Psi^\sharp + \nabla^\bot \Phi^\sharp$, and solving a Poisson equation. Thus, we have
\begin{equation} \norme{H^s}{u^{\sharp,2}} \leq C \norme{H^s}{v^\sharp}. \label{sharp} \end{equation}

\subsubsection{Proof of \propref{approx}, part 2 \\ Estimates on the approximate solution}

The important estimates are those on the first term $U^{1}$, or more precisely, those on the wave packet of unstable modes of the Euler equation
$$ u^{i,0}(t,x,y) = \int_{\mathbb{R}} \varphi^{1}(k)e^{ikx}e^{\lambda(k)t}\tilde{v}(k;y)~dk , $$
with $\varphi^{1}$ supported in $I=]-k_{0}-\eta,-k_{0}+\eta[\cup]k_{0}-\eta,k_{0}+\eta[$, where $\eta$ is small enough to have $I\subset \{k~|~\sigma(k)>0\}$, and to allow us to apply \thref{resth} in what follows. Also, $k_{0}$ must be the only critical point of $\sigma$ in $]k_{0}-\eta,k_{0}+\eta[$. By Parseval's equality, we have
$$ \norme{H^{s}}{u^{i,0}(t)}^{2} = \int_{\mathbb{R}} \varphi^{1}(k)^{2} \norme{H^{s}}{\tilde{v}(k)}^{2} e^{2\sigma(k)t}~dk $$
so, as $\pm k_{0}$ is the only nondegenerate maximum of $\sigma$ in each sub-interval of $I$, using the Laplace method, we get $\norme{H^{s}}{u^{i,0}(t)}^{2}\sim t^{-1/2}e^{2\sigma_{0}t}$ when $t\rightarrow +\infty$, so we have, given that $\norme{H^{s}}{u^{i,0}(t)}$ is bounded near $t=0$,
\begin{equation} C'\frac{e^{\sigma_{0}t}}{(1+t)^{1/4}} \leq \norme{H^{s}}{u^{i,0}(t)} \leq C\frac{e^{\sigma_{0}t}}{(1+t)^{1/4}} . \label{lapmet} \end{equation}

This allows us to estimate the other terms $u^{i,q}$ and $u^{b,q}$; this is done by induction. Fix $j$, and suppose that the internal layers $u^{i,q}$ with $0\leq q<j$ satisfy
\begin{equation} \norme{H^{s}}{u^{i,q}(t)} \leq c_{q}\frac{\exp\left[\sigma_{0}\left(1+\frac{q}{8n}\right)t\right]}{(1+t)^{1/4(1+q/8n)}} . \label{hruiq} \end{equation}
As $e^{\sigma_{0}\alpha t}(1+t)^{-\alpha}\stackrel{t\rightarrow+\infty}{\longrightarrow}+\infty$, we have
\begin{equation} \norme{H^{s}}{u^{i,q}(t)} \leq C_{j}\frac{\exp\left[\sigma_{0}\left(1+\frac{j}{8n}\right)t\right]}{(1+t)^{1/4(1+j/8n)}} . \label{uiq} \end{equation}
Then, with energy estimates on the heat equation with the Robin boundary condition induced by the Navier condition, we get
$$ \norme{H^{s}}{\frac{\bust(\sqrt\eps t)-\ust}{\eps^{1/8}}} \leq \frac{e^{C\sqrt{\eps}t}-1}{\eps^{1/8}}, $$
so $\eps^{-1/8}(\bust(\sqrt{\eps}t)-\ust)$ is bounded in $H^{s}$ uniformly for $t\leq \eps^{-1/32}$. Therefore, the remainder $r^{i,j}$ in the Euler equation $(E_{j})$ verifies (\ref{uiq}), as does $\partial_{t}r^{i,j}$, so $u^{i,j}$ also satisfies (\ref{uiq}) by \thref{resth}. Note that, because of our rescaling, we are studying the stability of $\bust(\sqrt\eps t)$.

Now we consider the boundary layer $u^{b,j}$. Suppose that $u^{b,q}$, $q<j$, are rapidly decreasing functions in the $Y$-variable and satisfy (\ref{hruiq}). By construction, $u^{b,j}_{2}$ and $p^{b,j}$ are rapidly decreasing functions in the $Y$-variable.
 Then, as $\eps^{-1/8}|\bust(\sqrt{\eps}t,\eps^{1/4}Y)|$ is also uniformly bounded in $\eps$ for $t\leq \eps^{-1/32}$, $Y\in\rplus$ and $\eps<\eps_{0}$ small enough, we get, by using the Green function of the heat equation with Neumann boundary conditions (see \cite{Sw}), that $u^{b,j}$ is a rapidly decreasing function of $Y$ satisfying
\begin{equation} \norme{H^{s}}{u^{b,j}} \leq C\sqrt{1+t}\frac{e^{\sigma_{0}(1+(j-1)/8n)t}}{(1+t)^{1/4[1+(j-1)/8n]}} \leq C_{j}\frac{\exp\left[\sigma_{0}\left(1+\frac{j}{8n}\right)t\right]}{(1+t)^{1/4(1+j/8n)}} . \label{ubiq} \end{equation}

Now, we estimate the main expansion terms, starting with $U^1$. Overall, after operating a change of variables in the boundary layer terms, we have
\begin{equation} \norme{H^{1}}{\eps^{n}U^{1}(t)} \leq C_{1}\sum_{j=0}^{8n-1} \frac{\eps^{n+j/8}}{(1+t)^{1/4(1+j/8n)}}\exp\left[\sigma_{0}\left(1+\frac{j}{8n}\right)t\right] . \label{eu1} \end{equation}
We will be interested in times $t=T^{\eps}_{0}-\tau$, with $T^{\eps}_{0}$ such that
$$ \frac{\eps^{n}e^{\sigma_{0}T^{\eps}_{0}}}{\sqrt{1+T^{\eps}_{0}}} = 1 . $$
The important manoeuvre here is to write $t=T^{\eps}_{0}-\tau$, so that
$$ \frac{\eps^{n}e^{\sigma_{0}t}}{\sqrt{1+t}} = e^{-\sigma_{0}\tau}\sqrt{\frac{1+T^{\eps}_{0}}{1+T^{\eps}_{0}-\tau}} := K^{\eps}_{0}(\tau)e^{-\sigma_{0}\tau} $$
with $K^{\eps}_{0}(\tau)\in ]1,2]$ for $\eps$ small enough, and with $0\leq\tau\leq\tau_{0}$ where $\tau_{0}$ does not depend on $\eps$. As $T^{\eps}_{0} \sim c\ln(\eps^{-1})$, we can choose $\eps<\eps_{0}$ so that $T^{\eps}_{0}-\tau_{0}>0$. $\tau_{0}$, and therefore $\eps_{0}$, will be chosen at the end of the proof in the next paragraph.

We show how this manoeuvre works in detail. Writing $t=T^{\eps}_{0}-\tau$, (\ref{eu1}) becomes
\begin{eqnarray*} \norme{H^{1}}{\eps^{n}U^{1}(T^{\eps}_{0}-\tau)} & \leq & (1+T^{\eps}_{0}-\tau)^{1/4}C_{1}K^{\eps}_{0}(\tau)e^{-\sigma_{0}\tau}\sum_{j=0}^{8n-1} (K^{\eps}_{0}(\tau)e^{-\sigma_{0}\tau})^{j/8n} \\ & \leq & (1+T^{\eps}_{0}-\tau)^{1/4}16nC_{1}K^{\eps}_{0}(\tau) e^{-\sigma_{0}\tau} . \end{eqnarray*}
Thus, by returning to the $t$-variable, for $t\leq T^{\eps}_{0}$, we get (\ref{est1}) for $U^1$:
$$ \norme{H^{1}}{U^{1}(t)} \leq (1+t)^{1/4}\frac{C_{1}e^{\sigma_{0}t}}{\sqrt{1+t}} . $$
This estimate means that $\norme{H^{1}}{U^{1}(t)}$ behaves like $\norme{H^{1}}{u^{i,0}(t)}$.
\newline

In \thref{theo3}, we want an estimate from below in $L^\infty$, so we will work a little more on $U^1$ in order to get a key under-estimate of a local $L^{2}$ norm, which will be (\ref{est2}). If $\alpha=\Im(\lambda'(k_{0}))$ and $\beta=-\sigma''(k_{0})>0$, then
$$ |u^{i,0}(t,x,y)| \sim e^{\sigma_{0}t} \left|\int_{I} \varphi^{1}(k)\tilde{v}(k,y)\exp(i(x+\alpha t)(k-k_{0})-\beta t (k-k_{0})^{2})~dk \right| . $$
But, writing $\kappa=k-k_{0}$, we have (by factorising)
$$ \int_{\mathbb{R}} \exp(i(x+\alpha t)\kappa-\beta t \kappa^{2})~d\kappa = \frac{c}{\sqrt{t}}\exp\left(\frac{-(x+\alpha t)^{2}}{4t}\right) $$
with $c$ complex, so there exists $C$ such that
$$ |u^{i,0}(t,x,y)| \geq \frac{Ce^{\sigma_{0}t}}{\sqrt{1+t}}\exp\left(\frac{-(x+\alpha t)^{2}}{4t}\right) . $$
Integrating $|u^{i,0}(t)|^{2}$ on the bounded domain
$$ \Omega_{A}(t) = \{(x,y)~|~y\leq A ~\mathrm{and}~|x+\alpha t|\leq A\sqrt{1+t}\} $$
for $A$ large enough, and using the fact that $u^{i,0}$ is the dominant term in $U^{1}$, we get (\ref{est2}):
$$ \norme{L^{2}(\Omega_{A}(t))}{U^{1}(t)} \geq C'_{1}(1+t)^{1/4}\frac{e^{\sigma_{0}t}}{\sqrt{1+t}} . $$
As the measure of $\Omega_{A}(t)$ is $A^{2}\sqrt{1+t}$, this gives us an under-estimate of the norm $\norme{L^{\infty}}{U^{1}(t)}$. In fact, $\norme{L^{\infty}}{U^{1}(t)}\sim C(1+t)^{-1/2}e^{\sigma_{0}t}$.
\newline

Now we specify the structure of $U^{j}$ to get $L^\infty$ and $H^1$ estimates: we write
$$ U^{j}(t,x,y) = \int_{I}\cdots\int_{I} V^{j}(k_{1},\cdots,k_{j};t,y) e^{ik_1 x}\cdots e^{ik_j x}~ dk_{1}\cdots dk_{j} $$
with $|V^{j}(k_{1},\cdots,k_{j};t,y)|\leq C_j\exp[(\sigma(k_{1})+\cdots+\sigma(k_{j}))t]$, which gives us
\begin{eqnarray*} \norme{L^{\infty}}{U^{j}(t)} & \leq & C_j\int_{I}\cdots\int_{I} \exp[(\sigma(k_{1})+\cdots+\sigma(k_{j}))t]~dk_{1}\cdots dk_{j} \\
 & \leq & C_j e^{j\sigma_{0}t} \prod_{m=1}^{j} \int_{\mathbb{R}}\exp(-\beta(k_m-k_{0})^{2}t)~dk_m \leq \frac{C_je^{j\sigma_{0}t}}{(1+t)^{j/2}} , \end{eqnarray*}
which is (\ref{estinf}), thanks to the second-order Taylor inequality
\begin{equation} \sigma(k) \leq \sigma_{0} - \beta(k-k_{0})^{2} . \label{tay} \end{equation}
For an $H^{1}$ estimate, we rewrite $U^{j}$ more precisely as
$$ U^{j}(t,x,y) = \int_{jk\in jI} \int_{k_{1}+\cdots+k_{j}=jk} V^{j}(k_{1},\cdots,k_{j};t,y)e^{ijkx}~dk_1\cdots dk_j $$
with $jI=I+\cdots+I$. As, for a given $k$, the value of $k_j$ is imposed by the other variables, Parseval's equality yields
$$ \norme{H^{1}}{U^{j}(t)}^{2} = \int_{jk\in jI} \norme{H^{1}}{\int_{k_{1}+\cdots+k_{j}=jk} V^{j}(k_{1},\cdots,k_{j};t)~dk_1\cdots dk_{j-1}}^{2}~dk := \int_{jk\in jI} {\cal N}(k,t)^{2}~dk . $$
Using (\ref{tay}), and noticing that
$$ \sum_{m=1}^{j} (k_{m}-k_{0})^{2} = j(k-k_{0})^{2} + \sum_{m=1}^{j} (k-k_{m})^{2} , $$
we have
$$ {\cal N}(t)^{2} \leq Ce^{2(j\sigma_{0}-j\beta(k-k_{0})^{2})t}\int_{\sum_{m=1}^j k_{m}=jk} \exp(-2\beta\sum_{m=1}^j (k_{m}-k)^{2}t)~dk_1\cdots dk_{j-1} . $$
Remembering that $k_j=jk-\sum_{m=1}^{j-1} k_m$, we integrate these gaussian functions, so
$$ \norme{H^{1}}{U^{j}(t)}^{2} \leq \int_{\mathbb{R}} \frac{Ce^{2j\sigma_{0}t-2j\beta(k-k_{0})^{2}t}}{t^{j-1}} ~dk . $$
Integrating this final gaussian function, and taking into account the boundedness of $\norme{H^{1}}{U^{j}(t)}$ near $t=0$, we get (\ref{est1}):
$$ \norme{H^{1}}{U^{j}(t)} \leq C_{j}(1+t)^{1/4}\frac{e^{j\sigma_{0}t}}{(1+t)^{j/2}} . $$

Finally, to end the proof of \propref{approx}, we need to prove that $\uap$ is indeed an approximate solution to the Navier-Stokes equation, in that the error $R$ satisfies estimate (\ref{Rest}). This error, $R=R^1 + R^\sharp$, contains two types of terms.
 Firstly, $R^1$ is made up of the terms $\eps^{(j+l)n}U^{j}\cdot \nabla U^{l}$ with $j+l>N$, in other words terms from the nonlinearity of the Navier-Stokes equation, so by using (\ref{451}), we have the estimate
$$ \norme{L^{2}}{R^1(t)} \leq C(1+t)^{1/4}\eps^{(N+1)n}\sum_{j=N+1}^{2N} \frac{e^{j\sigma_{0}t}}{(1+t)^{j/2}} . $$
Replacing $t$ by $T^{\eps}_{0}-\tau$, we get
\begin{eqnarray*} \norme{L^{2}}{R^1(T^{\eps}_{0}-\tau)} & \leq & C(1+T^{\eps}_{0}-\tau)^{1/4}\cdot K^{\eps}_0(\tau)^{2N}\cdot \sum_{j=N+1}^{2N} e^{-j\sigma_{0}\tau} \\
 & \leq & C(1+T^{\eps}_{0}-\tau)^{1/4}\cdot 2^{N-1}K^{\eps}_{0}(\tau)^{N+1}\cdot N e^{-j(N+1)\tau} , \end{eqnarray*}
and so, for $t\in[T^{\eps}_{0}-\tau_{0},T^{\eps}_{0}]$,
\begin{equation} \norme{L^2}{R^1(t)} \leq \frac{C\eps^{(N+1)n}}{(1+t)^{(N+1)/2}}(1+t)^{1/4}e^{(N+1)\sigma_{0}t}. \label{Rest2} \end{equation}
The constant $C$ is then adjusted so that (\ref{Rest2}) is true on $[0,T^\eps_0]$.

Secondly, in $R^\sharp$, we find terms involving the corrector term $\eps^{(N+1)n+l_0}u^\sharp$ and the error terms coming from the laplacian (the ones that do not appear in any of the equations $(E_q)$ or $(S_q)$) also appear. Specifically, we need estimates on the final four layers of $U^N$, so we set
\begin{eqnarray*} w^\sharp(t,x,y) & = & w^{i,\sharp}(t,x,y)+\eps^{1/4} w^{b,\sharp}(t,x,\eps^{-1/4}y) \\
 & = & \eps^{(N+1)n+l_0} \sum_{q=0}^{3} \left( \eps^{-q/8} u^{i,M-q} + \eps^{1/4-q/8} u^{b,M-q}(t,x,\eps^{-1/4}y) \right), \end{eqnarray*}
and we need to get estimates for $\norme{H^1}{w^\sharp}+\eps^{1/2}\norme{H^2}{w^\sharp}$. Indeed, given our construction of $u^{\sharp,1}$ and (\ref{sharp}), we have, due to the boundary-layer scale of $w^{b,\sharp}$,
$$ \norme{H^s}{\eps^{(N+1)n+l_0} u^\sharp} \leq C \eps^{3/8-s/4}\norme{H^s}{w^{b,\sharp}}. $$
We now note that, for any $l_0\geq 1$, the terms in $w^{i,\sharp}$ satisfy (\ref{uiq}) and those in $w^{b,\sharp}$ satisfy (\ref{ubiq}), both with $1+\frac{j}{8n}>N+1$, so we have that $w^\sharp$, and therefore $R^\sharp(t)$, satisfies (\ref{Rest2}) for $t\in [0,T^\eps_0]$, and therefore, we have (\ref{Rest}).

\subsubsection{Proof of \thref{theo3}: instability}

For this paragraph, we denote $\delta=\eps^{n}$. The order $N$ of the main expansion will be chosen in this section. Let $\ueps$ be an exact solution to the $\eqref{NS'}$ system with initial condition $\uap(0)$. We have $\norme{L^{2}}{\ueps(0)-\ust}\leq \eps^{n}$. The aim is to get a lower bound of $\norme{L^{\infty}}{\ueps-\bust}$, so we will look for a local $L^{2}$ lower bound.

Let us first estimate $\norme{L^{2}(\Omega_{A}(t))}{\uap-\bust}$: using (\ref{est1}) and (\ref{est2}), we have
$$ \norme{L^{2}(\Omega_{A}(t))}{\uap(t)-\bust(t)} \geq (1+t)^{1/4}\left(\frac{C_{0}\delta}{(1+t)^{1/2}}e^{\sigma_{0}t}-\sum_{j=2}^{N} \frac{C_{j}\delta^{j}}{(1+t)^{j/2}}e^{j\sigma_{0}t}\right) . $$
Writing $t=T^{\eps}_{0}-\tau$, we get
$$ \norme{L^{2}(\Omega_{A}(T^{\eps}_{0}-\tau))}{\uap-\bust} \geq (1+T^{\eps}_{0}-\tau)^{1/4}K^{\eps}_{0}(\tau)\left(C_{0}e^{-\sigma_{0}\tau} - C_{0}'e^{-2\sigma_{0}\tau}\right) $$
with $C_{0}'=N2^{N}\max_{1\leq j\leq N} C_{j}$. There exists $\tau_{1}>0$ such that
$$ C_{0}e^{-\sigma_{0}\tau}-C_{0}'e^{-2\sigma_{0}\tau} \geq \frac{C_{0}}{2} e^{-\sigma_{0}\tau} $$
for $\tau\geq \tau_{1}$. So, choosing $\tau_{0}>\tau_{1}$, we have, for $T^{\eps}_{0}-\tau_{0}\leq t\leq T^{\eps}_{0}-\tau_{1}$,
\begin{equation} \norme{L^{2}(\Omega_{A}(t))}{\uap(t)-\bust(t)} \geq (1+t)^{1/4}\frac{C_{0}\delta}{2(1+t)^{1/2}}e^{\sigma_{0}t} . \label{inst1} \end{equation}
Note that $\tau_{1}$ does not depend on $\delta$.
\newline

Now we estimate $\norme{L^{2}(\Omega)}{w(t)}$, where $w=\ueps-\uap$ and $q=\peps-\pap$. The pair $(w,q)$ solves the system
$$ \left\{ \begin{array}{rl} \partial_{t}w+w\cdot\nabla w-\sqrt{\eps}\Delta w+\nabla q & = -\uap\cdot\nabla w-w\cdot\nabla \uap +R \\
\div w & = 0 \\ (\partial_y w_1-2aw_1)(t,x,0) & = r(t,x,0) \\ w_2(t,x,0) & = 0 \end{array} \right. $$
with $w|_{t=0}=0$. Note that there is an error in the Navier boundary condition:
$$ r(t,x,0)=\eps^{(N+1)n+l_0}(\partial_y u^\sharp_1-2au^\sharp_1-2av^\sharp_1)(t,x,0). $$
We compute the $L^{2}$ dot-product of this equation and $w$; integrating by parts and using the Young inequality with a parameter $\rho_a$ depending on $a$ on the boundary term $\int_{\dOmega} r w_1$ (to absorb $\int_{\dOmega} |w_1|^2$), and the Young inequality again on the source term $\int_{\Omega} Rw$, we get the energy estimate
$$ \partial_{t}\norme{L^{2}}{w}^{2}+2\sqrt{\eps}\norme{L^{2}}{\nabla w}^{2}+2a\sqrt{\eps}\int_{\dOmega}|w_{1}|^{2} \hspace{150pt} $$
$$ \hspace{90pt} \leq (2\norme{L^{\infty}}{\nabla\uap}+\frac{1}{2})\norme{L^{2}}{w}^{2}+\frac{1}{2\rho_a}\sqrt{\eps}\int_{\dOmega} |r|^2 + \frac{1}{2}\norme{L^{2}}{R}^{2}. $$
With the trace theorem, we re-use the estimates on $u^\sharp$ and $v^\sharp$ to notice that the $H^1$ norm of the Navier error $r$ has the same behaviour as the $L^2$ norm of $R$, so there exists $C_a>0$ such that
$$ \partial_t\norme{L^2}{w}^2 + 2\sqrt{\eps}\left(\norme{L^2}{\nabla w}^2 + a\int_{\dOmega} |w_1|^2\right) \leq (2\norme{L^\infty}{\nabla \uap}+\frac{1}{2})\norme{L^2}{w}^2+C_a\norme{L^2}{R}^2. $$

An $L^{\infty}$ estimate on $\nabla\uap$ is crucial here: as each $U^{j}$ can be decomposed into two parts $U^{j}=U^{i,j}(t,x,y)+\eps^{1/4}U^{b,j}(t,x,Y)$, the first order derivatives of $U^{j}$ verify (\ref{estinf})-type estimates, since $\eps^{1/4}\nabla(U^{b,j}(t,x,\eps^{-1/4}y))=(\nabla U^{b,j})(t,x,Y)$, so
$$ \norme{L^{\infty}}{\nabla\uap(t)} \leq \norme{L^{\infty}}{\bust(t)}+\sum_{j=1}^{N} \frac{C_{j}\delta^{j}e^{j\sigma_{0}t}}{(1+t)^{j/2}} . $$
Using the $t=T^{\eps}_{0}-\tau$ manoeuvre, the sum on the right is smaller than $\frac{1}{4}$ for $t\leq T^{\eps}_{0}-\tau_{2}$, with $\tau_{2}\geq \tau_{1}$ independent of $\delta$.

Here we choose $N$, the order of the main expansion, to be such that, for $t\geq 0$, $2N\sigma_{0}>2(\norme{L^{\infty}}{\nabla\bust(t)}+1)$. With (\ref{Rest}) and this choice of $N$, the energy estimate becomes
$$ \partial_{t}\norme{L^{2}}{w(t)}^{2} \leq 2N\sigma_{0}\norme{L^{2}}{w(t)}^{2} + \frac{C_{R}\delta^{N+1}}{(1+t)^{(N+1)/2}}(1+t)^{1/4}e^{(N+1)\sigma_{0}t} . $$
Therefore, $\ueps(t)$ is an $L^{2}$ perturbation of $\bust(t)$, as we have the estimate
\begin{equation} \norme{L^{2}(\Omega)}{w(t)} \leq \frac{C_{0}'\delta^{N+1}}{(1+t)^{(N+1)/2}} (1+t)^{1/4} e^{(N+1)\sigma_{0}t} , \label{inst2} \end{equation}
by applying the following lemma, which is obtained using a variant of Gr\"onwall's lemma from \cite{PSV} and integration by parts:

\begin{lemma} Let $\varphi$ be a function such that
$$ \partial_{t}\varphi(t) \leq \lambda\varphi(t)+C\frac{e^{\mu t}}{(1+t)^{\alpha}} $$
for every $t\geq 0$, and for parameters $\mu>\lambda\geq 0$ and $\alpha>1$, then
$$ \varphi(t) \leq C'\frac{e^{\mu t}}{(1+t)^{\alpha}} $$
with $C'$ depending on $\varphi(0)$, $\lambda$, $\mu$ and $C$. \label{lemPSV} \end{lemma}

With (\ref{inst1}) and (\ref{inst2}), we can now conclude:
\begin{eqnarray*} \norme{L^{2}(\Omega_{A}(t))}{\ueps(t)-\bust(t)} & \geq & \norme{L^{2}(\Omega_{A}(t))}{\uap(t)-\bust(t)} - \norme{L^{2}(\Omega)}{w(t)} \\
 & \geq & (1+t)^{1/4}\left[\frac{C_{0}\delta e^{\sigma_{0}t}}{2(1+t)^{1/2}} - \frac{C_{0}'\delta^{N+1}e^{(N+1)\sigma_{0}t}}{(1+t)^{(N+1)/2}}\right] . \end{eqnarray*}
Writing $t=T^{\eps}_{0}-\tau$ for the last time, we have
$$ \frac{C_{0}\delta e^{\sigma_{0}t}}{2(1+t)^{1/2}} - \frac{C_{0}'\delta^{N+1}e^{(N+1)\sigma_{0}t}}{(1+t)^{(N+1)/2}} \geq \frac{C_{0}\delta e^{\sigma_{0}t}}{4(1+t)^{1/2}} $$
for $t\leq T^{\eps}_{0}-\tau_{3}$, with $\tau_{3}\geq \tau_{2}$ independent of $\delta$. Setting $\Teps=T^{\eps}_{0}-T$, with $\tau_{3}<T<\tau_{0}$ fixed (here we can finally fix $\tau_{0}$ and $\eps_{0}$), we have
$$ \norme{L^{2}(\Omega_{A}(\Teps))}{\ueps(\Teps)-\bust(\Teps)} \geq (1+T^{\eps})^{1/4}\frac{C_{0}}{4}e^{-\sigma_{0}T} := (1+T^{\eps})^{1/4}\deltaZ , $$
where $\deltaZ$ does not depend on $\eps$. Thus, $\norme{L^{\infty}}{\ueps(\Teps)-\bust(\Teps)} \geq C\deltaZ$.
\newline

Also, as $\dot{H}^{s}\hookrightarrow L^{\infty}$ for $s>1$, we have $\norme{\dot{H}^{s}}{\ueps(\Teps)-\bust(\Teps)}\geq C\deltaZ$. Remember that we are still using the fast variables, so, returning to the original scale of time and space (again without changing notation), we get (\ref{inst}) and
$$ \norme{\dot{H}^{s}}{\ueps(\Teps,x,y)-\bust(\Teps,\eps^{-1/2}y)} \geq C\deltaZ \eps^{(1-s)/2} \stackrel{\eps\rightarrow0}{\longrightarrow} +\infty . $$
Note that in the original scale of time, $\Teps = {\cal O}(n\sqrt{\eps}\ln(\eps^{-1}))$. $\square$
\newline

\preu{Remark:} we can also use the Sobolev embedding $\dot{H}^{s}\hookrightarrow L^{p^{*}}$ for $0<s<1$, with $p^{*}=\frac{2}{1-s}$, and the H\"older inequality to get a result similar to Grenier's in the Dirichlet case:
$$ \norme{\dot{H}^{s}(\Omega)}{\ueps(\Teps,x,y)-\bust\left(\Teps,\frac{y}{\sqrt{\eps}}\right)} \geq \deltaZ(s) \eps^{(1-s)/4} , $$
where $\deltaZ(s)$ is a constant depending on $s$, $\deltaZ$ in \thref{theo3} and $A$, the parameter in the bounded domain $\Omega_{A}(t)$. This gives a little information on the transition between stability in $L^{2}$ and instability in $L^{\infty}$.

\section{An example of unstable profile}

Examples of piecewise-linear flows that are linearly unstable for the Euler equation are given in \cite{DR} and \cite{Ge}, the latter stating that close-enough regularisations of these are smooth unstable shear flows. We wish to provide an explicit example of smooth linearly unstable profile, which is the object of \propref{exprop} below.
 From this example, we easily deduce the expression of a smooth unstable profile that fits \thref{theo3}, and we can also derive an unstable profile that tends to 0 exponentially. This final example fits the hypotheses of both of our main theorems in the periodic case, meaning that the convergence in $L^2$ and the instability in $L^\infty$ can be simultaneous.

\subsection{The model case: a hyperbolic tangent}

\begin{propo} For every $\delta>0$ and $\zeta\in\mathbb{R}$, the shear flow $(\sdel{u},0)$ with
$$ \sdel{u}(y) = \tanh(y-\delta)+\zeta $$
is a linearly unstable for the Euler equation. \label{exprop} \end{propo}

\preu{Proof:} writing $\lambda=-ikc$ for $\Re(\lambda)>0$ and $k\neq 0$, and taking the curl of the Euler equation linearised around $\sdel{u}$, we obtain a linear second-order differential equation for $\Psi$: the Rayleigh equation
$$ \rayl{c,k}: \hspace{10pt} (\sdel{u}-c)(\partial^{2}_{yy}-k^{2})\Psi - \sdel{u}''\Psi = 0 $$
with the conditions $\Psi(0)=0$ and $\lim_{y\rightarrow+\infty} \Psi(y)=0$. Our problem is now finding $c$ in the complex upper-half-plane, and real numbers $k$ such that $\rayl{c,k}$ has a solution in $H^{1}_{0}(\rplus)\cap H^{2}(\rplus)$.

It is well known that if the Rayleigh equation has an unstable solution, then $\sdel{u}$ must have an inflection point in $]0,+\infty[$ (Rayleigh's theorem, \cite{LR}), and $\sdel{u}$ has exactly one inflection point, $y_{0}=\delta$, with inflection value $u_{0}=\sdel{u}(y_{0})=\zeta$. Also, the function
$$ \sdel{K}(y) = \frac{-\sdel{u}''(y)}{\sdel{u}(y)-u_{0}} = \frac{-(\tanh(y-\delta))''}{\tanh(y-\delta)} = 2(1-\tanh(y-\delta)^{2}) $$
has a limit when $y\rightarrow y_{0}$, and is a positive continuous function on $\rplus$, vanishing at infinity. Thus, choosing $c=u_{0}$, we can divide $\rayl{c,k}$ by $\ust-u_{0}$ (which is not usually possible when $c$ is in the range of $\sdel{u}$), and we have a Sturm-Liouville problem
$$ -\Psi''(y) - \sdel{K}(y)\Psi(y) = -k^{2}\Psi(y) , $$
in which the square of the wave number $k$ intervenes as an eigenvalue of the operator $\sdel{S}=-\partial^{2}-\sdel{K}$. We shall therefore use the following result by Z.Lin (Theorem 1.5 in \cite{Lz}):

\begin{theo} \cite{Lz} Let $U$ be a ${\cal C}^{2}$ profile which has a limit $l$ as $y\rightarrow+\infty$, and an inflection point $y_{0}$ such that, writing $u_{0}=U(y_{0})$, the function $K(y)=\frac{-U''(y)}{U(y)-u_{0}}$ is a positive continuous function which goes to zero as $y\rightarrow+\infty$.
 We assume that $U$ takes the value $l$ only a finite number of times. If the lowest eigenvalue of the operator $S=-\partial^{2}-K$, defined on the Sobolev space $H^{1}_{0}(\rplus)\cap H^{2}(\rplus)$, is strictly negative, then, letting $-\mu^{2}$ be the lowest eigenvalue, for every $k\in(0,\mu)$, there exists $c(k)$ such that the Rayleigh equation $\rayl{c(k),k}$ has an unstable solution. \label{theo2} \end{theo}

The result is first shown on a finite interval; one proves that eigenfunctions of $S$ with negative eigenvalues are limits of unstable solutions, and this is done by using the Picard fixed point theorem. A compactness argument is then used to get the result on the half-line. We detail the proof of this theorem no further.
\newline

In the case of $\sdel{u}$, it remains to show that $\sdel{S}=-\partial^{2}-\sdel{K}$ has a negative eigenvalue. We hope to use a standard Sturm-Liouville argument, and notice that $\sdel{v}(y)=\tanh(y-\delta)$ is a solution to the equation
\begin{equation} -\varphi''(y)-\sdel{K}(y)\varphi(y) = 0 \label{SLd} \end{equation}
that has exactly one zero in $]0,+\infty[$. Since $\sdel{K}$ decreases exponentially to 0, the multiplication by $\sdel{K}$ is a compact perturbation of $-\partial^{2}$, so by Weyl's theorem (see \cite{RS}), the essential spectrum of $\sdel{S}$ is $\rplus$.
 But that means trouble, because the 0 in (\ref{SLd}) corresponds to the lowest point in the essential spectrum, and the one zero of $\sdel{v}$ means that $\sdel{S}$ has between 0 and 3 negative eigenvalues, according to the classical Sturm-Liouville analysis shown in \cite{DS}, chapter XIII.

So we must use a different tool to determine that $\sdel{S}$ has a strictly negative eigenvalue. Fix $\delta$ and let $\sdel{Q}$ be the quadratic form associated with $\sdel{S}$: for $u\in H^{1}_{0}(\rplus)$, we define
$$ \sdel{Q}(u) = \int_{0}^{+\infty} |u'(y)|^{2}-\sdel{K}(y)|u(y)|^{2}~dy . $$
Let $\seta{v}(y)=\tanh(y-\eta)$ for $\eta>0$, and $\chi$ be a [0,2]-supported smooth function, with $\chi(y)=1$ for $y\leq 1$. For any $n\in\mathbb{N}^{*}$, we define
$$ \seta{w}^{n}(y) = \left\{ \begin{array}{ll} 0 & \mathrm{if}~y\leq \eta \\ \seta{v}(y) & \mathrm{if}~y\in [\eta, \delta] \\
\seta{v}(y)\chi\left(\frac{y}{n}\right) & \mathrm{if}~y\geq\delta . \end{array} \right. $$
$\seta{w}^{n}$ is equal to $\seta{v}$ on $[\delta,n]$, continuous on $\rplus$, and is in $H^{1}_{0}(\rplus)$, and since $\seta{v}'$ and $\sqrt{\sdel{K}}\seta{v}$ are square-integrable on $\rplus$, we have
$$ \lim_{n\rightarrow+\infty} \sdel{Q}(\seta{w}^{n}) = \int_{\eta}^{+\infty} |\seta{v}'|^{2}-\sdel{K}|\seta{v}|^{2} := Q(\eta) . $$
By integrating by parts, we have $Q(\delta)=0$, and $Q$ is a differentiable function of $\eta$, so we look at $\partial_{\eta}Q(\delta)$. Defining $\seta{K}=-\seta{v}''/\seta{v}$, we have
$$ \partial_{\eta}Q(\eta) = \int_{\eta}^{+\infty} -2\seta{v}''\seta{v}'+2\sdel{K}\seta{v}\seta{v}'~dy - |\seta{v}'(\eta)|^2 = \int_{\delta}^{+\infty} 2\seta{v}\seta{v}'(\seta{K}+\sdel{K})~dy - 1 . $$
Thus, $\partial_\eta Q(\delta) = 8\int_\delta^{+\infty} \sdel{v}(y)\sdel{v}'(y)^2~dy -1$, and a change of variables $z=\sdel{v}(y+\delta)$ yields
$$ \partial_\eta Q(\delta) = 8\int_0^1 z(1-z^2) ~ dz - 1 = 1 > 0. $$
So there exists $\eta_{0}$ such that, for $\eta\in]\eta_{0},\delta[$, $Q(\eta)<0$, and therefore, for a given $\eta$ in that interval, there exists $n$ large enough such that $\seta{w}^{n}$ verifies $\sdel{Q}(\seta{w}^{n})<0$.

We have proved that the lowest point of the spectrum of $\sdel{S}$ is negative, and $\sdel{u}$ satisfies all the hypotheses of \thref{theo2}. \propref{exprop} is proved. $\square$
\newline

\preu{Remark:} $\sdel{S}$ has in fact exactly one negative eigenvalue. Indeed, remember that $\sdel{v}(y)=\tanh(y-\delta)$ is a bounded continuous solution of $-\varphi''-\sdel{K}\varphi=0$ that vanishes only at $y_{0}=\delta$, and that $\sdel{Q}(u)=\int_{0}^{+\infty} |u'|^{2}-\sdel{K}|u|^{2}$ is the quadratic form on $H^{1}_{0}(\rplus)$ associated with the operator $\sdel{S}$.
 We shall show that $\sdel{Q}\geq 0$ on $\sdel{F}=\{u\in H^{1}_{0}(\rplus)~|~ u(\delta)=0\}$. A function $u\in \sdel{F}$ can be written as $u=\sdel{v} w$, with $w\in H^{1}_{0}(\rplus)$. Let us take $w\in {\cal C}^{\infty}_{C}(\rplus)$, with $\mathrm{supp}(w)\subset [h,y_{0}-h]\cup[y_{0}+h,+\infty[$ for a certain $h>0$. Replacing $u$ by $\sdel{v} w$ in $Q(u)$, we get
$$ -u''-\sdel{K}u = -\sdel{v}''w-2\sdel{v}'w'-\sdel{v} w''-K\sdel{v} w = -2\sdel{v}'w'-\sdel{v} w'' $$
so, using $\sdel{Q}(u)=\int_{0}^{+\infty} \sdel{S}u\cdot u$, we have $\sdel{Q}(\sdel{v} w) = \int_{0}^{+\infty} -2\sdel{v}\sdel{v}' ww'-\sdel{v}^{2}ww''$. Integrating by parts, we get the factorisation $\sdel{Q}(\sdel{v} w) = \int_{0}^{+\infty} \sdel{v}^{2}w'^{2}$, \textit{i.e.}
\begin{equation} \sdel{Q}(u) = \int_{0}^{+\infty} \sdel{v}(y)^{2}\left|\left(\frac{u(y)}{\sdel{v}(y)}\right)'\right|^{2}~dy \hspace{10pt} \mathrm{for~all~} u\in \sdel{F} \label{21} \end{equation}
as ${\cal C}^{\infty}_{C}((0,\delta)\cup(\delta,+\infty[)$ is dense in $\sdel{F}$. We now have $\sdel{Q}\geq 0$ on $\sdel{F}$, which is a hyperplane of $H^{1}_{0}(\rplus)$, thus $\sdel{Q}$ is negative only on a subspace of dimension 1: $\sdel{S}$ has only one negative eigenvalue.
\newline

We can now give an example of unstable shear flow for \thref{theo3}: for a fixed positive $a$, it is one of the linearly unstable profiles for the Euler equation, $\sdel{u}(y)=\tanh(y-\delta)+\zeta$, with any $\delta>0$ and $\zeta$ to be chosen such that $\sdel{u}$ satisfies the rescaled Navier condition, $\frac{1}{2}\partial_{y}\sdel{u}(0) = a\sdel{u}(0)$. Setting $X=\tanh(\delta)$, we simply have
$$ \zeta = X+\frac{1}{2a}(1-X^{2}) > 0 . $$

\subsection{An $L^2$ example and simultaneous realisation of the results}

The idea is to ``truncate'' the hyperbolic tangent in a way that \thref{theo2} can still be applied. The starting point will therefore be a profile of the form
$$ u_\delta^\nu(y) = u_\delta(y) \xi(\nu y) , $$
where $\xi$ is ${\cal C}^\infty$, such that $\xi(y)=1$ for $y\leq 1$, decreasing for $y>1$ and goes to zero exponentially at infinity (compactly-supported will not do). For $\nu >0$ small enough, we will have $\xi(\nu y)=1$ for $y\in [0,2n_0]$, where $n_0$ is such that $Q_\delta(w_\eta^{n_0}) < 0$, as shown in the previous paragraph.
 Given that the Sturm-Liouville potential for the truncated profile will be equal to $K_\delta$ on the support of $w_\eta^{n_0}$, the quadratic forms for $w_\eta^{n_0}$ coincide and we will easily have that the lower bound on the spectrum of the corresponding Sturm-Liouville operator will be negative.

We must, however, be careful as to how we truncate the hyperbolic tangent. Indeed, to be able to apply Z.Lin's result, we must have that
$$ K_\delta^\nu = \frac{-(u_\delta^\nu)''}{u_\delta^\nu-\zeta} $$
is a positive (this is why a compactly-supported $\xi$ isn't convenient) and continuous function. This is not guaranteed because
$$ u_\delta^\nu - \zeta = \tanh(y-\delta)\xi(\nu y) - \zeta (1-\xi(\nu y)) $$
vanishes a second time, for $\nu y>1$. Given the flatness of the hyperbolic tangent when $y$ is large, it is reasonable to assume that it does not vanish a subsequent time, and that $u^\nu_\delta$ is decreasing for $\nu y \geq 5/4$. We must therefore construct $\xi$ so that this zero of the denominator of $K_\delta^\nu$ coincides with the second inflection point of $u_\delta^\nu$, which is generated by the truncation function $\xi$.
\newline

Set $\nu>0$ small enough so that $\xi(\nu y)=1$ on $[0,2n_0]$. We can assume that $n_0$ is large enough so that we have the inequality
$$ 1-\nu^2 < \tanh(y-\delta) < 1 $$
for $y\geq 2n_0$. This will allow us to get estimates on where the denominator will vanish. Indeed, we have
\begin{eqnarray} u_\delta^\nu(y)-\zeta & \leq & \xi(\nu y)-\zeta(1-\xi(\nu y)) \label{zeta1} \\
\mathrm{and}~ u_\delta^\nu(y) - \zeta & \geq & (1-\nu^2)\xi(\nu y) - \zeta(1-\xi(\nu y)) . \label{zeta2} \end{eqnarray}
The right-hand-sides of (\ref{zeta1}) and (\ref{zeta2}) are zero for a certain $y_0$ such that
$$ \frac{\zeta}{1+\zeta} < \xi(\nu y_0) < \frac{\zeta}{1+\zeta-\nu}. $$
Let us arbitrarily set $\xi(3/2) = \frac{\zeta}{1+\zeta}$; as $\xi$ is decreasing, the denominator therefore vanishes at $\nu y_0=\frac{3}{2} - h(\nu)$, with $h(\nu)>0$ going to zero as $\nu$ goes to zero.
 Indeed, impose that $\xi'$ has uniform-in-$\nu$ bounds on $[5/4,3/2]$ ($\zeta$ is a fixed number that does not depend on $\nu$, so this is reasonable) and that $\xi'<0$ on this interval. In this case, if $\nu$ is small enough so that $\frac{\zeta}{1+\zeta-\nu}<\xi(5/4)$, we get an estimate of $h(\nu)$ from first-order Taylor inequalities for $\xi$:
$$ |h(\nu)| \leq C(\xi(3/2-h)-\xi(3/2)) \leq C\zeta\left(\frac{1}{1+\zeta-\nu}-\frac{1}{1+\zeta}\right) \leq C\nu . $$

Then, by setting $\xi'(3/2-h)=-p$ for a fixed $p>0$ and by detailing the expression of $(u_\delta^\nu)''(y)$, we can get the value of $\xi''(3/2-h)$ that cancels out the numerator of $K_\delta^\nu$:
$$ \xi''(\nu y_0) = \frac{\xi(\nu y_0)}{\zeta} \left[ 2p(1-\tanh^2(y_0-\delta)) + 2\xi(\nu y_0) \tanh(y_0-\delta)(1-\tanh^2(y_0-\delta))\right] . $$
As $u_\delta^\nu$ is ${\cal C}^\infty$, we can then prescribe a third-derivative value that ensures that $K_\delta^\nu$ has a positive limit at $\nu y_0=3/2-h$ ($\xi^{(3)}(\nu y_0)<-2p\zeta$ for example). We have therefore constructed a function $\xi$ such that $K_\delta^\nu$ is defined, continuous and positive for $\nu y\leq 3/2$. We then set, for $z\geq 2$,
$$ \xi(z) = \frac{e^{-\gamma z}}{u_\delta (z/\nu)}, $$
so that $u_\delta^\nu(y) = e^{-\gamma\nu y}$ for $\nu y \geq 2$, thus
$$ K_\delta^\nu(y) = \frac{-\gamma^2 e^{-\gamma \nu y}}{e^{-\gamma\nu y} - \zeta} >0 $$
for $\nu y\geq 2$. We finally choose $\gamma$ large enough so that a ${\cal C}^\infty$ attachment to the first part of $\xi$ (for $\nu y\leq 3/2$) can be done while maintaining $\xi$ decreasing, and on which the convexity and the sign of $u_\delta^\nu-\zeta$ remain the same. As $\zeta$ is independent of $\nu$ and $\xi(3/2)$ is a fixed function of $\zeta$, $\gamma$ can be chosen independent of $\nu$.
\newline

Thus, for each $0<\nu<\nu_0$, we can construct a function $\xi^\nu$, element of a family $(\xi^\nu)_{\nu>0}$ which is uniformly bounded in ${\cal C}^3$ as $\nu\rightarrow 0$, so that $K_\delta^\nu$ is continuous and positive, equal to $K_\delta$ on $[0,2n_0]$ so that $Q_\delta^\nu(w_\eta^{n_0}) = Q_\delta(w_\eta^{n_0}) <0$, where, naturally,
$$ Q_\delta^\nu(u) = \int_0^{+\infty} |u'(y)|^2 - K_\delta^\nu(y) |u(y)|^2~dy $$
is the quadratic form associated with the operator $-\partial_y^2-K_\delta^\nu$ defined on $H^1_0(\rplus)\cap H^2(\rplus)$, proving that this operator has negative eigenvalues. We can therefore apply \thref{theo2} once again, and conclude with the following.
\begin{coro} For any $0<\nu<\nu_0$ and $\delta>0$, there exist $\zeta>0$ and a positive, non-increasing ${\cal C}^\infty$ function $\xi$ such that $\xi(y)= 1$ for $y\in [0,1]$ and $\xi(y)\stackrel{y\rightarrow+\infty}{\longrightarrow} 0$ at an exponential rate, such that $u_\delta^\nu(y) = (\tanh(y-\delta)+\zeta)\xi(\nu y)$ has two inflection points, and is an unstable profile for \thref{theo3}. \end{coro}

As it exponentially goes to zero at infinity, we have $u_\delta^\nu\in L^2(\rplus)$, and therefore $u_s(y)=(u_\delta^\nu(\eps^{-1/2} y),0)$ is in $L^2(\Omega')$, where $\Omega' = \mathbb{T}\times\rplus$.
 As both Theorems \ref{theo1} and \ref{theo3} are true in $\Omega'$ (see \cite{Ge} for the minor technical differences in the proof of the instability theorem), we get that not only are our two main results not contradictory, but, that in the periodic case, there are families of initial conditions that allow both to be satisfied simultaneously.
\vspace{12pt}

\noindent \textbf{\underline{Acknowledgements.}} I wish to express my warmest gratitude to my Ph.D supervisor, Fr\'ed\'eric Rousset, for the opportunity to work on this rich topic, and for his guidance and patience throughout the preparation of this article. 

\begin{small}
\bibliography{nsebib}
\bibliographystyle{abbrv}
\end{small}

\end{document}